\documentclass[12pt] {article}
\usepackage{amsmath,amsthm}
\usepackage{amssymb,latexsym}
\usepackage{amscd}
\title{Theory of valuations on manifolds, II.}
\date{}
\author{ Semyon Alesker \footnote{Partially supported by ISF grant 1369/04.}
\\  { \normalsize Department of Mathematics, Tel Aviv University, Ramat Aviv}
 \\  { \normalsize 69978 Tel Aviv,
Israel }
\\ {\normalsize e-mail: semyon@post.tau.ac.il}}

\def\eps{\varepsilon}
\def\alp{\alpha}
\def\ome{\omega}
\def\Ome{\Omega}
\def\lam{\lambda}
\def\Lam{\Lambda}

\def\to{\rightarrow}
\def\qed { Q.E.D. }
\def\pt{\partial}

\def\RR{\mathbb{R}}
\def\CC{\mathbb{C}}

\def\ZZ{\mathbb{Z}}

\def\PP{\mathbb{P}}

\def\an{ \oid {A}{}{} ( \CC^{n} ) }

\swapnumbers
\newtheorem{theorem}{Theorem}[subsection]
\newtheorem{corollary}[theorem]{Corollary}
\newtheorem{lemma}[theorem]{Lemma}
\newtheorem{proposition}[theorem]{Proposition}
\newtheorem{claim}[theorem]{Claim}
\theoremstyle{definition}
\newtheorem{example}[theorem]{Example}
\newtheorem{definition}[theorem]{Definition}
\newtheorem{remark}[theorem]{Remark}
\newtheorem{conjecture}[theorem]{Conjecture}
\theoremstyle{proposition-definition}
\newtheorem{proposition-definition}[theorem]{Proposition-Definition}

\def\cf{{\cal F}}

\def\ca{{\cal A}}  \def\cc{{\cal C}}
\def\cd{{\cal D}}  \def\cf{{\cal F}}
  
 \def\ck{{\cal K}} \def\cl{{\cal L}}
\def\cm{{\cal M}}  \def\co{{\cal O}}
\def\cp{{\cal P}} 
\def\car{{\cal R}}
\def\cs{{\cal S}} \def\ct{{\cal T}} \def\cu{{\cal U}}
\def\cv{{\cal V}} \def\cw{{\cal W}} 
 
\def\svv{SV(V)}

\def\qvv{QV (V)}

\makeatletter

\def\diagram{\m@th\leftwidth=\z@ \rightwidth=\z@ \topheight=\z@
\botheight=\z@ \setbox\@picbox\hbox\bgroup}

\def\enddiagram{\egroup\wd\@picbox\rightwidth\unitlength
\ht\@picbox\topheight\unitlength \dp\@picbox\botheight\unitlength
\hskip\leftwidth\unitlength\box\@picbox}

\def\bfig{\begin{diagram}}
\def\efig{\end{diagram}}
\newcount\wideness \newcount\leftwidth \newcount\rightwidth
\newcount\highness \newcount\topheight \newcount\botheight

\def\ratchet#1#2{\ifnum#1<#2 \global #1=#2 \fi}

\def\putbox(#1,#2)#3{%
\horsize{\wideness}{#3} \divide\wideness by 2
{\advance\wideness by #1 \ratchet{\rightwidth}{\wideness}}
{\advance\wideness by -#1 \ratchet{\leftwidth}{\wideness}}
\vertsize{\highness}{#3} \divide\highness by 2
{\advance\highness by #2 \ratchet{\topheight}{\highness}}
{\advance\highness by -#2 \ratchet{\botheight}{\highness}}
\put(#1,#2){\makebox(0,0){$#3$}}}

\def\putlbox(#1,#2)#3{%
\horsize{\wideness}{#3}
{\advance\wideness by #1 \ratchet{\rightwidth}{\wideness}}
{\ratchet{\leftwidth}{-#1}}
\vertsize{\highness}{#3} \divide\highness by 2
{\advance\highness by #2 \ratchet{\topheight}{\highness}}
{\advance\highness by -#2 \ratchet{\botheight}{\highness}}
\put(#1,#2){\makebox(0,0)[l]{$#3$}}}

\def\putrbox(#1,#2)#3{%
\horsize{\wideness}{#3}
{\ratchet{\rightwidth}{#1}}
{\advance\wideness by -#1 \ratchet{\leftwidth}{\wideness}}
\vertsize{\highness}{#3} \divide\highness by 2
{\advance\highness by #2 \ratchet{\topheight}{\highness}}
{\advance\highness by -#2 \ratchet{\botheight}{\highness}}
\put(#1,#2){\makebox(0,0)[r]{$#3$}}}

\def\adjust[#1]{} 

\newcount \coefa
\newcount \coefb
\newcount \coefc
\newcount\tempcounta
\newcount\tempcountb
\newcount\tempcountc
\newcount\tempcountd
\newcount\xext
\newcount\yext
\newcount\xoff
\newcount\yoff
\newcount\gap%
\newcount\arrowtypea
\newcount\arrowtypeb
\newcount\arrowtypec
\newcount\arrowtyped
\newcount\arrowtypee
\newcount\height
\newcount\width
\newcount\xpos
\newcount\ypos
\newcount\run
\newcount\rise
\newcount\arrowlength
\newcount\halflength
\newcount\arrowtype
\newdimen\tempdimen
\newdimen\xlen
\newdimen\ylen
\newsavebox{\tempboxa}%
\newsavebox{\tempboxb}%
\newsavebox{\tempboxc}%

\newdimen\w@dth

\def\setw@dth#1#2{\setbox\z@\hbox{\m@th$#1$}\w@dth=\wd\z@
\setbox\@ne\hbox{\m@th$#2$}\ifnum\w@dth<\wd\@ne \w@dth=\wd\@ne \fi
\advance\w@dth by 1.2em}

\def\t@^#1_#2{\allowbreak\def\n@one{#1}\def\n@two{#2}\mathrel
{\setw@dth{#1}{#2}
\mathop{\hbox to \w@dth{\rightarrowfill}}\limits
\ifx\n@one\empty\else ^{\box\z@}\fi
\ifx\n@two\empty\else _{\box\@ne}\fi}}
\def\t@@^#1{\@ifnextchar_{\t@^{#1}}{\t@^{#1}_{}}}
\def\to{\@ifnextchar^{\t@@}{\t@@^{}}}

\def\t@left^#1_#2{\def\n@one{#1}\def\n@two{#2}\mathrel{\setw@dth{#1}{#2}
\mathop{\hbox to \w@dth{\leftarrowfill}}\limits
\ifx\n@one\empty\else ^{\box\z@}\fi
\ifx\n@two\empty\else _{\box\@ne}\fi}}
\def\t@@left^#1{\@ifnextchar_{\t@left^{#1}}{\t@left^{#1}_{}}}
\def\toleft{\@ifnextchar^{\t@@left}{\t@@left^{}}}

\def\two@^#1_#2{\allowbreak
\def\n@one{#1}\def\n@two{#2}\mathrel{\setw@dth{#1}{#2}
\mathop{\vcenter{\lineskip\z@\baselineskip\z@
                 \hbox to \w@dth{\rightarrowfill}%
                 \hbox to \w@dth{\rightarrowfill}}%
       }\limits
\ifx\n@one\empty\else ^{\box\z@}\fi
\ifx\n@two\empty\else _{\box\@ne}\fi}}
\def\tw@@^#1{\@ifnextchar _{\two@^{#1}}{\two@^{#1}_{}}}
\def\two{\@ifnextchar ^{\tw@@}{\tw@@^{}}}

\def\tofr@^#1_#2{\def\n@one{#1}\def\n@two{#2}\mathrel{\setw@dth{#1}{#2}
\mathop{\vcenter{\hbox to \w@dth{\rightarrowfill}\kern-1.7ex
                 \hbox to \w@dth{\leftarrowfill}}%
       }\limits
\ifx\n@one\empty\else ^{\box\z@}\fi
\ifx\n@two\empty\else _{\box\@ne}\fi}}
\def\t@fr@^#1{\@ifnextchar_ {\tofr@^{#1}}{\tofr@^{#1}_{}}}
\def\tofro{\@ifnextchar^ {\t@fr@}{\t@fr@^{}}}

\def\mon{\mathop{\m@th\hbox to
      14.6\P@{\lasyb\char'51\hskip-2.1\P@$\arrext$\hss
$\mathord\rightarrow$}}\limits} 
\def\leftmono{\mathrel{\m@th\hbox to
14.6\P@{$\mathord\leftarrow$\hss$\arrext$\hskip-2.1\P@\lasyb\char'50%
}}\limits} 
\mathchardef\arrext="0200       

\def\settypes(#1,#2,#3){\arrowtypea#1 \arrowtypeb#2 \arrowtypec#3}
\def\settoheight#1#2{\setbox\@tempboxa\hbox{#2}#1\ht\@tempboxa\relax}%
\def\settodepth#1#2{\setbox\@tempboxa\hbox{#2}#1\dp\@tempboxa\relax}%
\def\settokens`#1`#2`#3`#4`{%
     \def\tokena{#1}\def\tokenb{#2}\def\tokenc{#3}\def\tokend{#4}}
\def\setsqparms[#1`#2`#3`#4;#5`#6]{%
\arrowtypea #1
\arrowtypeb #2
\arrowtypec #3
\arrowtyped #4
\width #5
\height #6
}
\def\setpos(#1,#2){\xpos=#1 \ypos#2}

\def\settriparms[#1`#2`#3;#4]{\settripairparms[#1`#2`#3`1`1;#4]}%

\def\settripairparms[#1`#2`#3`#4`#5;#6]{%
\arrowtypea #1
\arrowtypeb #2
\arrowtypec #3
\arrowtyped #4
\arrowtypee #5
\width #6
\height #6
}

\def\resetparms{\settripairparms[1`1`1`1`1;500]\width 500}

\resetparms

\def\mvector(#1,#2)#3{
\put(0,0){\vector(#1,#2){#3}}%
\put(0,0){\vector(#1,#2){26}}%
}
\def\evector(#1,#2)#3{{
\arrowlength #3
\put(0,0){\vector(#1,#2){\arrowlength}}%
\advance \arrowlength by-30
\put(0,0){\vector(#1,#2){\arrowlength}}%
}}

\def\horsize#1#2{%
\settowidth{\tempdimen}{$#2$}%
#1=\tempdimen
\divide #1 by\unitlength
}

\def\vertsize#1#2{%
\settoheight{\tempdimen}{$#2$}%
#1=\tempdimen
\settodepth{\tempdimen}{$#2$}%
\advance #1 by\tempdimen
\divide #1 by\unitlength
}

\def\putvector(#1,#2)(#3,#4)#5#6{{%
\ifnum3<\arrowtype
\putdashvector(#1,#2)(#3,#4)#5\arrowtype
\else
\ifnum\arrowtype<-3
\putdashvector(#1,#2)(#3,#4)#5\arrowtype
\else
\xpos=#1
\ypos=#2
\run=#3
\rise=#4
\arrowlength=#5
\ifnum \arrowtype<0
    \ifnum \run=0
        \advance \ypos by-\arrowlength
    \else
        \tempcounta \arrowlength
        \multiply \tempcounta by\rise
        \divide \tempcounta by\run
        \ifnum\run>0
            \advance \xpos by\arrowlength
            \advance \ypos by\tempcounta
        \else
            \advance \xpos by-\arrowlength
            \advance \ypos by-\tempcounta
        \fi
    \fi
    \multiply \arrowtype by-1
    \multiply \rise by-1
    \multiply \run by-1
\fi
\ifcase \arrowtype
\or \put(\xpos,\ypos){\vector(\run,\rise){\arrowlength}}%
\or \put(\xpos,\ypos){\mvector(\run,\rise)\arrowlength}%
\or \put(\xpos,\ypos){\evector(\run,\rise){\arrowlength}}%
\fi\fi\fi
}}

\def\putsplitvector(#1,#2)#3#4{
\xpos #1
\ypos #2
\arrowtype #4
\halflength #3
\arrowlength #3
\gap 140
\advance \halflength by-\gap
\divide \halflength by2
\ifnum\arrowtype>0
   \ifcase \arrowtype
   \or \put(\xpos,\ypos){\line(0,-1){\halflength}}%
       \advance\ypos by-\halflength
       \advance\ypos by-\gap
       \put(\xpos,\ypos){\vector(0,-1){\halflength}}%
   \or \put(\xpos,\ypos){\line(0,-1)\halflength}%
       \put(\xpos,\ypos){\vector(0,-1)3}%
       \advance\ypos by-\halflength
       \advance\ypos by-\gap
       \put(\xpos,\ypos){\vector(0,-1){\halflength}}%
   \or \put(\xpos,\ypos){\line(0,-1)\halflength}%
       \advance\ypos by-\halflength
       \advance\ypos by-\gap
       \put(\xpos,\ypos){\evector(0,-1){\halflength}}%
   \fi
\else \arrowtype=-\arrowtype
   \ifcase\arrowtype
   \or \advance \ypos by-\arrowlength
       \put(\xpos,\ypos){\line(0,1){\halflength}}%
       \advance\ypos by\halflength
       \advance\ypos by\gap
       \put(\xpos,\ypos){\vector(0,1){\halflength}}%
   \or \advance \ypos by-\arrowlength
       \put(\xpos,\ypos){\line(0,1)\halflength}%
       \put(\xpos,\ypos){\vector(0,1)3}%
       \advance\ypos by\halflength
       \advance\ypos by\gap
       \put(\xpos,\ypos){\vector(0,1){\halflength}}%
   \or \advance \ypos by-\arrowlength
       \put(\xpos,\ypos){\line(0,1)\halflength}%
       \advance\ypos by\halflength
       \advance\ypos by\gap
       \put(\xpos,\ypos){\evector(0,1){\halflength}}%
   \fi
\fi
}

\def\putmorphism(#1)(#2,#3)[#4`#5`#6]#7#8#9{{%
\run #2
\rise #3
\ifnum\rise=0
  \puthmorphism(#1)[#4`#5`#6]{#7}{#8}#9%
\else\ifnum\run=0
  \putvmorphism(#1)[#4`#5`#6]{#7}{#8}#9%
\else
\setpos(#1)%
\arrowlength #7
\arrowtype #8
\ifnum\run=0
\else\ifnum\rise=0
\else
\ifnum\run>0
    \coefa=1
\else
   \coefa=-1
\fi
\ifnum\arrowtype>0
   \coefb=0
   \coefc=-1
\else
   \coefb=\coefa
   \coefc=1
   \arrowtype=-\arrowtype
\fi
\width=2
\multiply \width by\run
\divide \width by\rise
\ifnum \width<0  \width=-\width\fi
\advance\width by60
\if l#9 \width=-\width\fi
\putbox(\xpos,\ypos){#4}
{\multiply \coefa by\arrowlength
\advance\xpos by\coefa
\multiply \coefa by\rise
\divide \coefa by\run
\advance \ypos by\coefa
\putbox(\xpos,\ypos){#5} }%
{\multiply \coefa by\arrowlength
\divide \coefa by2
\advance \xpos by\coefa
\advance \xpos by\width
\multiply \coefa by\rise
\divide \coefa by\run
\advance \ypos by\coefa
\if l#9%
   \putrbox(\xpos,\ypos){#6}%
\else\if r#9%
   \putlbox(\xpos,\ypos){#6}%
\fi\fi }%
{\multiply \rise by-\coefc
\multiply \run by-\coefc
\multiply \coefb by\arrowlength
\advance \xpos by\coefb
\multiply \coefb by\rise
\divide \coefb by\run
\advance \ypos by\coefb
\multiply \coefc by70
\advance \ypos by\coefc
\multiply \coefc by\run
\divide \coefc by\rise
\advance \xpos by\coefc
\multiply \coefa by140
\multiply \coefa by\run
\divide \coefa by\rise
\advance \arrowlength by\coefa
\ifcase\arrowtype
\or \put(\xpos,\ypos){\vector(\run,\rise){\arrowlength}}%
\or \put(\xpos,\ypos){\mvector(\run,\rise){\arrowlength}}%
\or \put(\xpos,\ypos){\evector(\run,\rise){\arrowlength}}%
\fi}\fi\fi\fi\fi}}

\newcount\numbdashes \newcount\lengthdash \newcount\increment

\def\howmanydashes{
\numbdashes=\arrowlength \lengthdash=40
\divide\numbdashes by \lengthdash
\lengthdash=\arrowlength
\divide\lengthdash by \numbdashes
\increment=\lengthdash
\multiply\lengthdash by 3
\divide\lengthdash by 5
}

\def\putdashvector(#1)(#2,#3)#4#5{%
\ifnum#3=0 \putdashhvector(#1){#4}#5
\else
\ifnum#2=0
\putdashvvector(#1){#4}#5\fi\fi}

\def\putdashhvector(#1,#2)#3#4{{%
\arrowlength=#3 \howmanydashes
\multiput(#1,#2)(\increment,0){\numbdashes}%
{\vrule height .4pt width \lengthdash\unitlength}
\arrowtype=#4 \xpos=#1
\ifnum\arrowtype<0 \advance\arrowtype by 7 \fi
\ifcase\arrowtype
\or \advance\xpos by 10
    \put(\xpos,#2){\vector(-1,0){\lengthdash}}
    \advance\xpos by 40
    \put(\xpos,#2){\vector(-1,0){\lengthdash}}
\or \advance \xpos by 10
    \put(\xpos,#2){\vector(-1,0){\lengthdash}}
    \advance\xpos by  \arrowlength
    \advance\xpos by  -50
    \put(\xpos,#2){\vector(-1,0){\lengthdash}}
\or \advance\xpos by 10
    \put(\xpos,#2){\vector(-1,0){\lengthdash}}
\or \advance\xpos by \arrowlength
    \advance\xpos by -\lengthdash
    \put(\xpos,#2){\vector(1,0){\lengthdash}}
\or {\advance\xpos by 10
    \put(\xpos,#2){\vector(1,0){\lengthdash}}}
    \advance\xpos by \arrowlength
    \advance\xpos by -\lengthdash
    \put(\xpos,#2){\vector(1,0){\lengthdash}}
\or \advance\xpos by \arrowlength
    \advance\xpos by -\lengthdash
    \put(\xpos,#2){\vector(1,0){\lengthdash}}
    \advance\xpos by -40
    \put(\xpos,#2){\vector(1,0){\lengthdash}}
   \fi
}}

\def\putdashvvector(#1,#2)#3#4{{%
\arrowlength=#3 \howmanydashes
\ypos=#2 \advance\ypos by -\arrowlength
\multiput(#1,#2)(0,\increment){\numbdashes}%
    {\vrule width .4pt height \lengthdash\unitlength}
\arrowtype=#4 \ypos=#2
\ifnum\arrowtype<0 \advance\arrowtype by 7 \fi
\ifcase\arrowtype
\or \advance\ypos by \arrowlength \advance\ypos by -40
    \put(#1,\ypos){\vector(0,1){\lengthdash}}
    \advance\ypos by -40
    \put(#1,\ypos){\vector(0,1){\lengthdash}}
\or \advance\ypos by 10
    \put(#1,\ypos){\vector(0,1){\lengthdash}}
    \advance\ypos by \arrowlength \advance\ypos by -40
    \put(#1,\ypos){\vector(0,1){\lengthdash}}
\or \advance\ypos by \arrowlength \advance\ypos by -40
    \put(#1,\ypos){\vector(0,1){\lengthdash}}
\or \advance\ypos by 10
    \put(#1,\ypos){\vector(0,-1){\lengthdash}}
\or \advance\ypos by 10
    \put(#1,\ypos){\vector(0,-1){\lengthdash}}
    \advance\ypos by \arrowlength \advance\ypos by -40
    \put(#1,\ypos){\vector(0,-1){\lengthdash}}
\or \advance\ypos by 10
    \put(#1,\ypos){\vector(0,-1){\lengthdash}}
    \advance\ypos by 40
    \put(#1,\ypos){\vector(0,-1){\lengthdash}}
\fi
}}

\def\puthmorphism(#1,#2)[#3`#4`#5]#6#7#8{{%
\xpos #1
\ypos #2
\width #6
\arrowlength #6
\arrowtype=#7
\putbox(\xpos,\ypos){#3\vphantom{#4}}%
{\advance \xpos by\arrowlength
\putbox(\xpos,\ypos){\vphantom{#3}#4}}%
\horsize{\tempcounta}{#3}%
\horsize{\tempcountb}{#4}%
\divide \tempcounta by2
\divide \tempcountb by2
\advance \tempcounta by30
\advance \tempcountb by30
\advance \xpos by\tempcounta
\advance \arrowlength by-\tempcounta
\advance \arrowlength by-\tempcountb
\putvector(\xpos,\ypos)(1,0)\arrowlength\arrowtype
\divide \arrowlength by2
\advance \xpos by\arrowlength
\vertsize{\tempcounta}{#5}%
\divide\tempcounta by2
\advance \tempcounta by20
\if a#8 %
   \advance \ypos by\tempcounta
   \putbox(\xpos,\ypos){#5}%
\else
   \advance \ypos by-\tempcounta
   \putbox(\xpos,\ypos){#5}%
\fi}}

\def\putvmorphism(#1,#2)[#3`#4`#5]#6#7#8{{%
\xpos #1
\ypos #2
\arrowlength #6
\arrowtype #7
\settowidth{\xlen}{$#5$}%
\putbox(\xpos,\ypos){#3}%
{\advance \ypos by-\arrowlength
\putbox(\xpos,\ypos){#4}}%
{\advance\arrowlength by-140
\advance \ypos by-70
\ifdim\xlen>0pt
   \if m#8%
      \putsplitvector(\xpos,\ypos)\arrowlength\arrowtype
   \else
   \putvector(\xpos,\ypos)(0,-1)\arrowlength\arrowtype
   \fi
\else
   \putvector(\xpos,\ypos)(0,-1)\arrowlength\arrowtype
\fi}%
\ifdim\xlen>0pt
   \divide \arrowlength by2
   \advance\ypos by-\arrowlength
   \if l#8%
      \advance \xpos by-40
      \putrbox(\xpos,\ypos){#5}%
   \else\if r#8%
      \advance \xpos by40
      \putlbox(\xpos,\ypos){#5}%
   \else
      \putbox(\xpos,\ypos){#5}%
   \fi\fi
\fi
}}

\def\putsquarep<#1>(#2)[#3;#4`#5`#6`#7]{{%
\setsqparms[#1]%
\setpos(#2)%
\settokens`#3`%
\puthmorphism(\xpos,\ypos)[\tokenc`\tokend`{#7}]{\width}{\arrowtyped}b%
\advance\ypos by \height
\puthmorphism(\xpos,\ypos)[\tokena`\tokenb`{#4}]{\width}{\arrowtypea}a%
\putvmorphism(\xpos,\ypos)[``{#5}]{\height}{\arrowtypeb}l%
\advance\xpos by \width
\putvmorphism(\xpos,\ypos)[``{#6}]{\height}{\arrowtypec}r%
}}

\def\putsquare{\@ifnextchar <{\putsquarep}{\putsquarep%
   <\arrowtypea`\arrowtypeb`\arrowtypec`\arrowtyped;\width`\height>}}
\def\square{\@ifnextchar< {\squarep}{\squarep
   <\arrowtypea`\arrowtypeb`\arrowtypec`\arrowtyped;\width`\height>}}
\def\squarep<#1>[#2`#3`#4`#5;#6`#7`#8`#9]{{
\setsqparms[#1]
\diagram
\putsquarep<\arrowtypea`\arrowtypeb`\arrowtypec`
\arrowtyped;\width`\height>
(0,0)[#2`#3`#4`{#5};#6`#7`#8`{#9}]
\enddiagram
}}                                                 
\def\putptrianglep<#1>(#2,#3)[#4`#5`#6;#7`#8`#9]{{%
\settriparms[#1]%
\xpos=#2 \ypos=#3
\advance\ypos by \height
\puthmorphism(\xpos,\ypos)[#4`#5`{#7}]{\height}{\arrowtypea}a%
\putvmorphism(\xpos,\ypos)[`#6`{#8}]{\height}{\arrowtypeb}l%
\advance\xpos by\height
\putmorphism(\xpos,\ypos)(-1,-1)[``{#9}]{\height}{\arrowtypec}r%
}}

\def\putptriangle{\@ifnextchar <{\putptrianglep}{\putptrianglep
   <\arrowtypea`\arrowtypeb`\arrowtypec;\height>}}
\def\ptriangle{\@ifnextchar <{\ptrianglep}{\ptrianglep
   <\arrowtypea`\arrowtypeb`\arrowtypec;\height>}}
\def\ptrianglep<#1>[#2`#3`#4;#5`#6`#7]{{
\settriparms[#1]
\diagram
\putptrianglep<\arrowtypea`\arrowtypeb`
\arrowtypec;\height>
(0,0)[#2`#3`#4;#5`#6`{#7}]
\enddiagram
}}                                            

\def\putqtrianglep<#1>(#2,#3)[#4`#5`#6;#7`#8`#9]{{%
\settriparms[#1]%
\xpos=#2 \ypos=#3
\advance\ypos by\height
\puthmorphism(\xpos,\ypos)[#4`#5`{#7}]{\height}{\arrowtypea}a%
\putmorphism(\xpos,\ypos)(1,-1)[``{#8}]{\height}{\arrowtypeb}l%
\advance\xpos by\height
\putvmorphism(\xpos,\ypos)[`#6`{#9}]{\height}{\arrowtypec}r%
}}

\def\putqtriangle{\@ifnextchar <{\putqtrianglep}{\putqtrianglep
   <\arrowtypea`\arrowtypeb`\arrowtypec;\height>}}
\def\qtriangle{\@ifnextchar <{\qtrianglep}{\qtrianglep
   <\arrowtypea`\arrowtypeb`\arrowtypec;\height>}}
\def\qtrianglep<#1>[#2`#3`#4;#5`#6`#7]{{
\settriparms[#1]
\width=\height                                
\diagram
\putqtrianglep<\arrowtypea`\arrowtypeb`
\arrowtypec;\height>
(0,0)[#2`#3`#4;#5`#6`{#7}]
\enddiagram
}}

\def\putdtrianglep<#1>(#2,#3)[#4`#5`#6;#7`#8`#9]{{%
\settriparms[#1]%
\xpos=#2 \ypos=#3
\puthmorphism(\xpos,\ypos)[#5`#6`{#9}]{\height}{\arrowtypec}b%
\advance\xpos by \height \advance\ypos by\height
\putmorphism(\xpos,\ypos)(-1,-1)[``{#7}]{\height}{\arrowtypea}l%
\putvmorphism(\xpos,\ypos)[#4``{#8}]{\height}{\arrowtypeb}r%
}}

\def\putdtriangle{\@ifnextchar <{\putdtrianglep}{\putdtrianglep
   <\arrowtypea`\arrowtypeb`\arrowtypec;\height>}}
\def\dtriangle{\@ifnextchar <{\dtrianglep}{\dtrianglep
   <\arrowtypea`\arrowtypeb`\arrowtypec;\height>}}
\def\dtrianglep<#1>[#2`#3`#4;#5`#6`#7]{{
\settriparms[#1]
\width=\height                                
\diagram
\putdtrianglep<\arrowtypea`\arrowtypeb`
\arrowtypec;\height>
(0,0)[#2`#3`#4;#5`#6`{#7}]
\enddiagram
}}

\def\putbtrianglep<#1>(#2,#3)[#4`#5`#6;#7`#8`#9]{{%
\settriparms[#1]%
\xpos=#2 \ypos=#3
\puthmorphism(\xpos,\ypos)[#5`#6`{#9}]{\height}{\arrowtypec}b%
\advance\ypos by\height
\putmorphism(\xpos,\ypos)(1,-1)[``{#8}]{\height}{\arrowtypeb}r%
\putvmorphism(\xpos,\ypos)[#4``{#7}]{\height}{\arrowtypea}l%
}}

\def\putbtriangle{\@ifnextchar <{\putbtrianglep}{\putbtrianglep
   <\arrowtypea`\arrowtypeb`\arrowtypec;\height>}}
\def\btriangle{\@ifnextchar <{\btrianglep}{\btrianglep
   <\arrowtypea`\arrowtypeb`\arrowtypec;\height>}}
\def\btrianglep<#1>[#2`#3`#4;#5`#6`#7]{{
\settriparms[#1]
\width=\height                               
\diagram
\putbtrianglep<\arrowtypea`\arrowtypeb`
\arrowtypec;\height>
(0,0)[#2`#3`#4;#5`#6`{#7}]
\enddiagram
}}

\def\putAtrianglep<#1>(#2,#3)[#4`#5`#6;#7`#8`#9]{{%
\settriparms[#1]%
\xpos=#2 \ypos=#3
{\multiply \height by2
\puthmorphism(\xpos,\ypos)[#5`#6`{#9}]{\height}{\arrowtypec}b}%
\advance\xpos by\height \advance\ypos by\height
\putmorphism(\xpos,\ypos)(-1,-1)[#4``{#7}]{\height}{\arrowtypea}l%
\putmorphism(\xpos,\ypos)(1,-1)[``{#8}]{\height}{\arrowtypeb}r%
}}

\def\putAtriangle{\@ifnextchar <{\putAtrianglep}{\putAtrianglep
   <\arrowtypea`\arrowtypeb`\arrowtypec;\height>}}
\def\Atriangle{\@ifnextchar <{\Atrianglep}{\Atrianglep
   <\arrowtypea`\arrowtypeb`\arrowtypec;\height>}}
\def\Atrianglep<#1>[#2`#3`#4;#5`#6`#7]{{
\settriparms[#1]
\width=\height                                     
\diagram
\putAtrianglep<\arrowtypea`\arrowtypeb`
\arrowtypec;\height>
(0,0)[#2`#3`#4;#5`#6`{#7}]
\enddiagram
}}

\def\putAtrianglepairp<#1>(#2)[#3;#4`#5`#6`#7`#8]{{%
\settripairparms[#1]%
\setpos(#2)%
\settokens`#3`%
\puthmorphism(\xpos,\ypos)[\tokenb`\tokenc`{#7}]{\height}{\arrowtyped}b%
\advance\xpos by\height
\puthmorphism(\xpos,\ypos)[\phantom{\tokenc}`\tokend`{#8}]%
{\height}{\arrowtypee}b%
\advance\ypos by\height
\putmorphism(\xpos,\ypos)(-1,-1)[\tokena``{#4}]{\height}{\arrowtypea}l%
\putvmorphism(\xpos,\ypos)[``{#5}]{\height}{\arrowtypeb}m%
\putmorphism(\xpos,\ypos)(1,-1)[``{#6}]{\height}{\arrowtypec}r%
}}

\def\putAtrianglepair{\@ifnextchar <{\putAtrianglepairp}{\putAtrianglepairp%
   <\arrowtypea`\arrowtypeb`\arrowtypec`\arrowtyped`\arrowtypee;\height>}}
\def\Atrianglepair{\@ifnextchar <{\Atrianglepairp}{\Atrianglepairp%
   <\arrowtypea`\arrowtypeb`\arrowtypec`\arrowtyped`\arrowtypee;\height>}}

\def\Atrianglepairp<#1>[#2;#3`#4`#5`#6`#7]{{
\settripairparms[#1]
\settokens`#2`
\width=\height                                
\diagram
\putAtrianglepairp                            
<\arrowtypea`\arrowtypeb`\arrowtypec`
\arrowtyped`\arrowtypee;\height>
(0,0)[{#2};#3`#4`#5`#6`{#7}]
\enddiagram
}}

\def\putVtrianglep<#1>(#2,#3)[#4`#5`#6;#7`#8`#9]{{%
\settriparms[#1]%
\xpos=#2 \ypos=#3
\advance\ypos by\height
{\multiply\height by2
\puthmorphism(\xpos,\ypos)[#4`#5`{#7}]{\height}{\arrowtypea}a}%
\putmorphism(\xpos,\ypos)(1,-1)[`#6`{#8}]{\height}{\arrowtypeb}l%
\advance\xpos by\height
\advance\xpos by\height
\putmorphism(\xpos,\ypos)(-1,-1)[``{#9}]{\height}{\arrowtypec}r%
}}

\def\putVtriangle{\@ifnextchar <{\putVtrianglep}{\putVtrianglep
   <\arrowtypea`\arrowtypeb`\arrowtypec;\height>}}
\def\Vtriangle{\@ifnextchar <{\Vtrianglep}{\Vtrianglep
   <\arrowtypea`\arrowtypeb`\arrowtypec;\height>}}
\def\Vtrianglep<#1>[#2`#3`#4;#5`#6`#7]{{
\settriparms[#1]
\width=\height                                 
\diagram
\putVtrianglep<\arrowtypea`\arrowtypeb`
\arrowtypec;\height>
(0,0)[#2`#3`#4;#5`#6`{#7}]
\enddiagram
}}

\def\putVtrianglepairp<#1>(#2)[#3;#4`#5`#6`#7`#8]{{
\settripairparms[#1]%
\setpos(#2)%
\settokens`#3`%
\advance\ypos by\height
\putmorphism(\xpos,\ypos)(1,-1)[`\tokend`{#6}]{\height}{\arrowtypec}l%
\puthmorphism(\xpos,\ypos)[\tokena`\tokenb`{#4}]{\height}{\arrowtypea}a%
\advance\xpos by\height
\puthmorphism(\xpos,\ypos)[\phantom{\tokenb}`\tokenc`{#5}]%
{\height}{\arrowtypeb}a%
\putvmorphism(\xpos,\ypos)[``{#7}]{\height}{\arrowtyped}m%
\advance\xpos by\height
\putmorphism(\xpos,\ypos)(-1,-1)[``{#8}]{\height}{\arrowtypee}r%
}}

\def\putVtrianglepair{\@ifnextchar <{\putVtrianglepairp}{\putVtrianglepairp%
    <\arrowtypea`\arrowtypeb`\arrowtypec`\arrowtyped`\arrowtypee;\height>}}
\def\Vtrianglepair{\@ifnextchar <{\Vtrianglepairp}{\Vtrianglepairp%
    <\arrowtypea`\arrowtypeb`\arrowtypec`\arrowtyped`\arrowtypee;\height>}}
\def\Vtrianglepairp<#1>[#2;#3`#4`#5`#6`#7]{{
\settripairparms[#1]
\settokens`#2`
\diagram
\putVtrianglepairp                             
<\arrowtypea`\arrowtypeb`\arrowtypec`
\arrowtyped`\arrowtypee;\height>
(0,0)[{#2};#3`#4`#5`#6`{#7}]
\enddiagram
}}

\def\putCtrianglep<#1>(#2,#3)[#4`#5`#6;#7`#8`#9]{{%
\settriparms[#1]%
\xpos=#2 \ypos=#3
\advance\ypos by\height
\putmorphism(\xpos,\ypos)(1,-1)[``{#9}]{\height}{\arrowtypec}l%
\advance\xpos by\height
\advance\ypos by\height
\putmorphism(\xpos,\ypos)(-1,-1)[#4`#5`{#7}]{\height}{\arrowtypea}l%
{\multiply\height by 2
\putvmorphism(\xpos,\ypos)[`#6`{#8}]{\height}{\arrowtypeb}r}%
}}

\def\putCtriangle{\@ifnextchar <{\putCtrianglep}{\putCtrianglep
    <\arrowtypea`\arrowtypeb`\arrowtypec;\height>}}
\def\Ctriangle{\@ifnextchar <{\Ctrianglep}{\Ctrianglep
    <\arrowtypea`\arrowtypeb`\arrowtypec;\height>}}
\def\Ctrianglep<#1>[#2`#3`#4;#5`#6`#7]{{
\settriparms[#1]
\width=\height                               
\diagram
\putCtrianglep<\arrowtypea`\arrowtypeb`
\arrowtypec;\height>
(0,0)[#2`#3`#4;#5`#6`{#7}]
\enddiagram
}}                                           
\def\putDtrianglep<#1>(#2,#3)[#4`#5`#6;#7`#8`#9]{{%
\settriparms[#1]%
\xpos=#2 \ypos=#3
\advance\xpos by\height \advance\ypos by\height
\putmorphism(\xpos,\ypos)(-1,-1)[``{#9}]{\height}{\arrowtypec}r%
\advance\xpos by-\height \advance\ypos by\height
\putmorphism(\xpos,\ypos)(1,-1)[`#5`{#8}]{\height}{\arrowtypeb}r%
{\multiply\height by 2
\putvmorphism(\xpos,\ypos)[#4`#6`{#7}]{\height}{\arrowtypea}l}%
}}

\def\putDtriangle{\@ifnextchar <{\putDtrianglep}{\putDtrianglep
    <\arrowtypea`\arrowtypeb`\arrowtypec;\height>}}
\def\Dtriangle{\@ifnextchar <{\Dtrianglep}{\Dtrianglep
   <\arrowtypea`\arrowtypeb`\arrowtypec;\height>}}
\def\Dtrianglep<#1>[#2`#3`#4;#5`#6`#7]{{
\settriparms[#1]
\width=\height                              
\diagram
\putDtrianglep<\arrowtypea`\arrowtypeb`
\arrowtypec;\height>
(0,0)[#2`#3`#4;#5`#6`{#7}]
\enddiagram
}}                                          
\def\setrecparms[#1`#2]{\width=#1 \height=#2}%

\def\recursep<#1`#2>[#3;#4`#5`#6`#7`#8]{{\m@th
\width=#1 \height=#2
\settokens`#3`
\settowidth{\tempdimen}{$\tokena$}
\ifdim\tempdimen=0pt
  \savebox{\tempboxa}{\hbox{$\tokenb$}}%
  \savebox{\tempboxb}{\hbox{$\tokend$}}%
  \savebox{\tempboxc}{\hbox{$#6$}}%
\else
  \savebox{\tempboxa}{\hbox{$\hbox{$\tokena$}\times\hbox{$\tokenb$}$}}%
  \savebox{\tempboxb}{\hbox{$\hbox{$\tokena$}\times\hbox{$\tokend$}$}}%
  \savebox{\tempboxc}{\hbox{$\hbox{$\tokena$}\times\hbox{$#6$}$}}%
\fi
\ypos=\height
\divide\ypos by 2
\xpos=\ypos
\advance\xpos by \width
\bfig
\putCtrianglep<-1`1`1;\ypos>(0,0)[`\tokenc`;#5`#6`{#7}]%
\puthmorphism(\ypos,0)[\tokend`\usebox{\tempboxb}`{#8}]{\width}{-1}b%
\puthmorphism(\ypos,\height)[\tokenb`\usebox{\tempboxa}`{#4}]{\width}{-1}a%
\advance\ypos by \width
\putvmorphism(\ypos,\height)[``\usebox{\tempboxc}]{\height}1r%
\efig
}}

\def\recurse{\@ifnextchar <{\recursep}{\recursep<\width`\height>}}

\def\puttwohmorphisms(#1,#2)[#3`#4;#5`#6]#7#8#9{{%
\puthmorphism(#1,#2)[#3`#4`]{#7}0a
\ypos=#2
\advance\ypos by 20
\puthmorphism(#1,\ypos)[\phantom{#3}`\phantom{#4}`#5]{#7}{#8}a
\advance\ypos by -40
\puthmorphism(#1,\ypos)[\phantom{#3}`\phantom{#4}`#6]{#7}{#9}b
}}

\def\puttwovmorphisms(#1,#2)[#3`#4;#5`#6]#7#8#9{{%
\putvmorphism(#1,#2)[#3`#4`]{#7}0a
\xpos=#1
\advance\xpos by -20
\putvmorphism(\xpos,#2)[\phantom{#3}`\phantom{#4}`#5]{#7}{#8}l
\advance\xpos by 40
\putvmorphism(\xpos,#2)[\phantom{#3}`\phantom{#4}`#6]{#7}{#9}r
}}

\def\puthcoequalizer(#1)[#2`#3`#4;#5`#6`#7]#8#9{{%
\setpos(#1)%
\puttwohmorphisms(\xpos,\ypos)[#2`#3;#5`#6]{#8}11%
\advance\xpos by #8
\puthmorphism(\xpos,\ypos)[\phantom{#3}`#4`#7]{#8}1{#9}
}}

\def\putvcoequalizer(#1)[#2`#3`#4;#5`#6`#7]#8#9{{%
\setpos(#1)%
\puttwovmorphisms(\xpos,\ypos)[#2`#3;#5`#6]{#8}11%
\advance\ypos by -#8
\putvmorphism(\xpos,\ypos)[\phantom{#3}`#4`#7]{#8}1{#9}
}}

\def\putthreehmorphisms(#1)[#2`#3;#4`#5`#6]#7(#8)#9{{%
\setpos(#1) \settypes(#8)
\if a#9 %
     \vertsize{\tempcounta}{#5}%
     \vertsize{\tempcountb}{#6}%
     \ifnum \tempcounta<\tempcountb \tempcounta=\tempcountb \fi
\else
     \vertsize{\tempcounta}{#4}%
     \vertsize{\tempcountb}{#5}%
     \ifnum \tempcounta<\tempcountb \tempcounta=\tempcountb \fi
\fi
\advance \tempcounta by 60
\puthmorphism(\xpos,\ypos)[#2`#3`#5]{#7}{\arrowtypeb}{#9}
\advance\ypos by \tempcounta
\puthmorphism(\xpos,\ypos)[\phantom{#2}`\phantom{#3}`#4]{#7}{\arrowtypea}{#9}
\advance\ypos by -\tempcounta \advance\ypos by -\tempcounta
\puthmorphism(\xpos,\ypos)[\phantom{#2}`\phantom{#3}`#6]{#7}{\arrowtypec}{#9}
}}

\def\setarrowtoks[#1`#2`#3`#4`#5`#6]{%
\def\toka{#1}
\def\tokb{#2}
\def\tokc{#3}
\def\tokd{#4}
\def\toke{#5}
\def\tokf{#6}
}
\def\hex{\@ifnextchar <{\hexp}{\hexp<1000`400>}}
\def\hexp<#1`#2>[#3`#4`#5`#6`#7`#8;#9]{%
\setarrowtoks[#9]
\yext=#2 \advance \yext by #2
\xext=#1 \advance\xext by \yext
\bfig
\putCtriangle<-1`0`1;#2>(0,0)[`#5`;\tokb``\tokd]
\xext=#1 \yext=#2 \advance \yext by #2
\putsquare<1`0`0`1;\xext`\yext>(#2,0)[#3`#4`#7`#8;\toka```\tokf]
\advance \xext by #2
\putDtriangle<0`1`-1;#2>(\xext,0)[`#6`;`\tokc`\toke]
\efig
}
\makeatother

\def\inj{\hookrightarrow}
\def\surj{\twoheadrightarrow}
\begin{document}
\maketitle \setcounter{section}{-1}
\begin{abstract}
This article is the second part in the series of articles where we
are developing theory of valuations on manifolds. Roughly speaking
valuations could be thought as finitely additive measures on a
class of nice subsets of a manifold which satisfy some additional
assumptions.

The goal of this article is to introduce a notion of a smooth
valuation on an arbitrary smooth manifold and establish some of
the basic properties of it.
\end{abstract}
\tableofcontents
\section{Introduction.}\setcounter{subsection}{1}
This article is the second part in the series of articles where we
are developing theory of valuations on manifolds. Roughly speaking
valuations could be thought as finitely additive measures on a
class of nice subsets of a manifold which satisfy some additional
assumptions.

The goal of this article is to introduce the notion of a smooth
valuation on an arbitrary smooth manifold and establish some of
the basic properties of it. Let us describe this notion with
several oversimplifications referring for the details to the main
text.

Let $X$ be a smooth manifold of dimension $n$. Let us denote by
$\cp(X)$ the family of simple differentiable subpolyhedra of $X$
(see Subsection \ref{poly}). $\cp(X)$ serves as a natural class of
"nice" sets. For any set $P\in \cp(X)$ one defines a
characteristic cycle $CC(P)$ which is a closed cycle of dimension
$n$ in the cotangent bundle $T^*X$ (Definition
\ref{define-char-c}). (Note that if $P$ is a smooth submanifold of
$X$ then $CC(P)$ coincides with the conormal bundle of $P$.) A
smooth valuation $\phi$ is a complex valued finitely additive
functional (measure) on $\cp(X)$ which satisfies some additional
properties. One of the main such properties is continuity of
$\phi$ with respect to convergence in the sense of currents of the
characteristic cycles of subsets from $\cp(X)$. Most of the other
properties were introduced essentially for technical reasons and
their necessity is not very clear for the moment.

\begin{remark}
The class $\cp(X)$ is not closed neither under finite unions nor
under finite intersections. Thus the notion of a finitely additive
functional on $\cp(X)$ should be defined more formally. This is
done in Subsections \ref{fa} and \ref{sheaf1} using the notion of
a subdivision of a differentiable polyhedron.
\end{remark}

Thus we get the space $V^\infty(X)$ of smooth valuations on $X$.
It is a Fr\'echet space. The group of diffeomorphisms acts
continuously on $V^\infty(X)$. It is important to notice that if
$X$ is an affine space then the subspace of translation invariant
elements from $V^\infty(X)$ coincides with the space $Val^{sm}(X)$
introduced and studied by the author in \cite{alesker-poly}; the
last space is a dense subspace of the space $Val(X)$ of continuous
translation invariant valuations on {\itshape convex} subsets of
$X$ which is the classical object. For the classical theory of
valuations we refer to the surveys McMullen-Schneider
\cite{mcmullen-schneider} and McMullen \cite{mcmullen-survey}.

Next, the notion of smooth valuation is a local notion. More
precisely for any open subset $\cu\subset X$ the correspondence
$\cu\mapsto V^\infty (\cu)$ is a sheaf on $X$ (when the
restriction maps are obvious). This sheaf is denoted by
$\cv^\infty_X$. Thus $V^\infty(X)$ is equal to the space of global
sections $\Gamma(X,\cv^\infty_X)=\cv^\infty_X(X)$.

The sheaf $\cv^\infty_X$ has a canonical filtration by subsheaves
of vector spaces
$$\cv^\infty_X=\cw_0\supset \cw_1\supset \dots\supset \cw_n$$
where $n=\dim X$. $\cw_n$ coincides with the sheaf of smooth
densities (measures) on $X$. For any open subset $\cu\subset X$
and any $i=0,1,\dots,n$, $\cw_i(\cu)$ is a closed subspace of
$\cv^\infty_X(\cu)$.

It turns out that the associated graded sheaf
$gr_\cw\cv^\infty_X:=\bigoplus_{i=0}^n\cw_i/\cw_{i+1}$ admits a
simple description in terms of translation invariant valuations.
To state it let us denote by $Val(TX)$ the (infinite dimensional)
vector bundle over $X$ such that its fiber over a point $x\in X$
is equal to the space $Val^{sm}(T_xX)$ of smooth translation
invariant valuations of the tangent space $T_xX$. By McMullen's
Theorem \ref{val-3} the space $Val^{sm}(T_xX)$ has natural grading
by the degree of homogeneity which must be an integer between 0
and $n$. Thus $Val(TX)$ is a graded vector bundle. Let us denote
by $\underline{Val}(TX)$ the sheaf $\cu\mapsto
C^\infty(\cu,Val(TX))$ where the last space denotes the space of
infinitely smooth sections of $Val(TX)$ over $\cu$. The next
result is Corollary \ref{filt-7}.
\begin{theorem}
There exists a canonical isomorphism of graded sheaves
$$gr_\cw\cv^\infty_X\simeq \underline{Val}(TX).$$
\end{theorem}
This theorem provides a description of smooth valuations since
translation invariant valuations are studied much better.
Proposition \ref{filt-9} gives yet another description of smooth
valuations in terms of integration with respect to the
characteristic cycle. Combined with Lemma \ref{sheaf-sm-1} it says
the following.
\begin{theorem}
Let $\phi$ be a smooth valuation on $X$. Then there exists a
section $\eta\in \tilde C^\infty(T^*X,\Omega^n\otimes p^*o)$ such
that for any $P\in \cp(X)$ one has
$$\phi(P)=\int_{CC(P)}\eta$$
where $p:T^*X\to X$ is the canonical projection, $\Ome^n$ denotes
the vector bundle of $n$-forms on $T^*X$, $o$ denotes the
orientation bundle on $X$, and $\tilde
C^\infty(T^*X,\Omega^n\otimes p^*o)$ denotes the space of
infinitely smooth sections of the bundle $\Omega^n\otimes p^*o$
such that the restriction of the projection $p$ to the support of
this section is proper.

Conversely any expression of the above form is a smooth valuation.
\end{theorem}

The sheaf $\cv^\infty_X$ has yet another interesting structure
which we call the Euler-Verdier involution and denote by $\sigma$.
This is a non-trivial automorphism of sheaf $\sigma\colon
\cv^\infty_X \tilde\to \cv^\infty_X$. The next result is Theorem
\ref{eu-ve-th}.
\begin{theorem}
(i) The Euler-Verdier involution $\sigma$ preserves the filtration
$\cw_\bullet$.

(ii) The induced involution on $gr_\cw \cv_X^\infty \simeq
\underline{Val}_\bullet(TX)$ comes from the involution on the
bundle $Val(TX)$ defined as $\phi\mapsto [K\mapsto (-1)^{deg \phi}
\phi(-K)]$ for any $\phi\in Val(T_xX)$ for any $x\in X$, and where
$deg \phi$ is the degree of homogeneity of $\phi$.
\end{theorem}

Thus the sheaf $\cv_X^\infty$ of smooth valuations decomposes
under the action of the Euler-Verdier involution into two
subsheaves $\cv_X^{\infty, +}$ and $\cv_X^{\infty,-}$
corresponding to eigenvalues 1 and -1 of $\sigma$ respectively.
Thus
$$\cv_X^\infty=\cv_X^{\infty, +}\oplus \cv_X^{\infty, -}.$$

The article is organized as follows. Section \ref{background} is a
background and does not contain new results. In Subsection
\ref{rep} we remind very basic facts from representation theory,
in Subsection \ref{sheaf} we remind some basic facts from the
sheaf theory, and in Subsection \ref{valuations} facts from the
theory of valuations.

In Section \ref{pomeva} we discuss the notion and the properties
of differentiable polyhedra, discuss the notion of a finitely
additive measure on them, and finally in Subsection \ref{sheaf-sm}
we introduce the main object of this article, the notion of a
smooth valuation on a manifold.

In Section \ref{further} we study further general properties of
smooth valuations. In Subsection \ref{further-1} we introduce and
study the filtration $\cw_\bullet$ on smooth valuations; in
Proposition \ref{filt-9} we obtain a description of smooth
valuations in terms of the integration with respect to the
characteristic cycle. In Subsection \ref{topology} we define the
natural structure of Fr\'echet space on the space of smooth
valuations. In Subsection \ref{eu-ve} we introduce the
Euler-Verdier involution on smooth valuations.

{\bf Acknowledgements.}  I am grateful to J. Bernstein to numerous
very useful discussions. I express my gratitude to J. Fu for very
fruitful conversations, and in particular for his explanations of
the construction of valuations using the integration over the
normal (characteristic) cycle. I thank V.D. Milman for his
interest to this work, and D. Hug and M. Ludwig for useful
discussions. I thank the referee for careful reading the first
version of the article.

\section{Background}\label{background}\setcounter{subsection}{0}
In  Subsection \ref{rep} we remind some very basic definitions and
facts from representation theory. In Subsection \ref{sheaf} we
remind basic facts from the sheaf theory. In Subsection
\ref{valuations} we remind some facts from the valuation theory.
This section does not contain new results.

\subsection{Some representation theory.}
\label{rep} \setcounter{theorem}{0}
\begin{definition}\label{rep-1} Let $\rho$ be a continuous representation of a Lie group $G$ in a
Fr\'echet space $F$. A vector $\xi \in F$ is called $G$-smooth if
the map $g\mapsto \rho(g)\xi$ is an infinitely differentiable map
from $G$ to $F$.
\end{definition}
It is well known (see e.g. \cite{wallach}, Section 1.6) that the
subset $F^{sm}$ of smooth vectors is a $G$-invariant linear
subspace dense in $F$. Moreover it has a natural topology of a
Fr\'echet space (which is stronger than that induced from $F$),
and the representation of $G$ in $F^{sm}$ is continuous. Moreover
all vectors in $F^{sm}$ are $G$-smooth.
\subsection{Sheaf theory.}\label{sheaf} The definitions of this
subsection are taken from Godement's book \cite{godement}.

 Let $X$ be a topological space. Let
$\Phi$ be a family of closed subsets of $X$.
\begin{definition}[\cite{godement}, Section 3.2]\label{sh-1}
The family $\Phi$ is called paracompactifiable if

(1) any $S\in \Phi$ is closed and paracompact;

(2) $\Phi$ is closed under finite unions;

(3) any closed subset of any $S\in \Phi$ also belongs to $\Phi$;

(4) any $S\in \Phi$ has a neighborhood belonging to $\Phi$.
\end{definition}
From now on we will always assume that $\Phi$ is a
paracompactifiable family of subsets of $X$.

\begin{example}\label{sh-2}
(1) If $X$ is a locally compact paracompact topological space then
the family of all closed subsets is paracompactifiable.

(2) If $X$ is a locally compact paracompact topological space then
the family of all compact subsets is paracompactifiable.
\end{example}

From now on we will always assume that $X$ is locally compact and
paracompact.
\begin{definition}[\cite{godement}, Section 3.5]\label{sh-3}
(1) A sheaf $\cf$ on $X$ is called $\Phi$-soft if for any $S',\,
S''\in \Phi$ with $S'\supset S''$ the restriction map $\cf(S')\to
\cf(S'')$ is surjective.

(2) A sheaf $\cf$ on $X$ is called soft if it is $\Phi$-soft where
$\Phi$ is the family of all closed subsets of $X$.
\end{definition}

For a sheaf $\cf$ on $X$ let us denote by $\Gamma_\Phi(\cf)$ the
set of global sections of $\cf$ with the support in $\Phi$. The
functor $\cf\mapsto \Gamma_\Phi(\cf)$ is left exact on the
category of sheaves of abelian groups. Denote as usual by
$H^i_\Phi(X,\cf)$ its right derived functor.
\begin{theorem}[\cite{godement}, Theorem 3.5.4]\label{sh-4}
Let $0\to \cl^0\to \cl^1\to \dots$ be an exact sequence of
$\Phi$-soft sheaves of abelian groups. Then the following sequence
is exact:
$$0\to \Gamma_\Phi(\cl^0)\to \Gamma_\Phi(\cl^1)\to \dots.$$
\end{theorem}
\begin{theorem}[\cite{godement}, Theorem 4.4.3]\label{sh-5}
If $\cf$ is a $\Phi$-soft sheaf then $H^i_\Phi(X,\cf)=0 $ for all
$i>0$.
\end{theorem}
\begin{theorem}[\cite{godement}, Theorem 3.7.1]\label{sh-6}
Let $\ca$ be a sheaf of unital rings on $X$. If $\ca$ is
$\Phi$-soft then any $\ca$-module is $\Phi$-soft.
\end{theorem}

\begin{theorem}[\cite{godement}, Theorem 3.7.2]\label{sh-7}
Let $\ca$ be a sheaf of unital rings over a paracompact space $X$.
Then $\ca$ is soft if and only if any point of $X$ has a
neighborhood $U$ such that for any disjoint closed subsets
$S,T\subset U$ there exists a section of $\ca$ over $U$ which is
equal to 1 on $S$ and to 0 on $T$.
\end{theorem}

\begin{definition}[\cite{godement}, Section 3.7]\label{sh-8}
Let $\cl$ be a sheaf of abelian groups on $X$. $\cl$ is called
fine (resp. $\Phi$-fine) if the sheaf of rings ${\cal  H} om_
{\ZZ} (\cl,\cl)$ is soft (resp. $\cl |_S$ is fine for all $S\in
\Phi$).
\end{definition}

Any $\Phi$-fine sheaf is $\Phi$-soft (by Theorem \ref{sh-6}).
\begin{lemma}\label{sh-9}
Let $\ca$ be a sheaf of unital rings on a locally compact
paracompact space $X$. Then $\ca$ is $\Phi$-fine if and only if it
is $\Phi$-soft.
\end{lemma}
{\bf Proof.} By the remark before this lemma it remains to prove
that if $\ca$ is $\Phi$-soft then it is $\Phi$-fine. Restricting
to $S\in \Phi$, it suffices to prove that if $\ca$ is soft than it
is fine. Note that a sheaf $\cl$ of abelian groups is fine if and
only if for any disjoint closed subsets $A$ and $B$ of $X$ there
exists a morphism $\cl\to \cl$ inducing the identity map over a
neighborhood of $A$ and the zero map over a neighborhood of $B$.
Thus in the case of a sheaf of rings $\ca$, the last condition is
satisfied provided $\ca$ has a section over $X$ which is equal to
one in a neighborhood of $A$ and is equal to zero in a
neighborhood of $B$. The last condition is equivalent to the fact
that $\ca$ is soft. \qed

\begin{theorem}[\cite{godement}, Theorem 3.7.3]\label{sh-10}
If $\cl$ is a $\Phi$-fine sheaf of abelian groups then for any
sheaf $\cm$ of abelian groups $\cl\otimes_{\ZZ} \cm$ is
$\Phi$-fine (and hence $\Phi$-soft).
\end{theorem}
\begin{example}[\cite{godement}, Section 3.7]\label{sh-11}
Let $X$ be a smooth paracompact manifold. Let $\co_X$ denote the
sheaf of $C^\infty$-functions on $X$. Then $\co_X$ is fine and
hence soft. Hence any  $\co_X$-module $\cm$ is fine and soft. It
follows that $H^i(X,\cm)=H^i_c(X,\cm)=0$ for all $i>0$.
\end{example}

For the definition of operations $f_*,f^*,f_!$ on sheaves see e.g.
the book \cite{kashiwara-schapira} (in their notation $f^*$ is
denoted by $f^{-1}$).

\begin{lemma}\label{sh-12}
Let $X$ be a locally compact paracompact topological space. Let
$Z\subset X$ be a closed subset of $X$. Consider the imbeddings
$$Z\overset{i}{\hookrightarrow}X \overset{j}{\hookleftarrow}
U:=X\backslash Z.$$ Let $\cf$ be a sheaf on $X$.

1) If $H^1(U,j_!j^*\cf)=0$ then any section of $\cf$ over $Z$
extends to a section over $X$.

2) If any section of $\cf$ over $Z$ extends to a section over $X$
and $H^1(X,\cf)=0$ then $H^1(U,j_!j^*\cf)=0$.

3) Let $\ca$ be a soft sheaf of unital rings on $X$. Then for any
$\ca|_U$-module $\cm$ one has
$$H^i(X,j_!\cm)=0 \mbox{ for all } i>0.$$
\end{lemma}
{\bf Proof.} 1) We have  an exact sequence of sheaves
$$0\to j_!j^*\cf \to \cf \to i_*i^*\cf\to 0$$
(see e.g. the exact sequence in Proposition 2.3.6(v) of
\cite{kashiwara-schapira} combined with Propositions 2.3.6(iv) and
2.5.4(ii) of \cite{kashiwara-schapira}). Hence the following
sequence is exact
$$\Gamma(X,\cf)\to \Gamma(X,i_*i^*\cf)\to H^1(X,j_!j^*\cf)=0.$$
Since $\Gamma(X,i_*i^*\cf)=\Gamma(Z,i^*\cf)$ the result follows.

2) From the same exact sequence and our assumptions we obtain an
exact sequence
$$\Gamma(X,\cf)\overset{\mbox{onto}}{\twoheadrightarrow} \Gamma(Z,i^*\cf)\to
H^1(X,j_!j^*\cf)\to H^1(X,\cf)=0.$$ This implies the statement.

3) Indeed $j_!\cm$ is an $\ca$-module. Hence it is acyclic by
Theorems \ref{sh-6} and \ref{sh-5}. \qed

\begin{proposition}\label{sh-13}
Let $X$ be a smooth manifold. Let $\co_X$ denote the sheaf of
$C^\infty$-functions on $X$.  Let $\cv$ be a sheaf on $X$ which
admits a finite filtration by subsheaves
$$\cv=\cw_0\supset \cw_1 \supset \dots \supset \cw_N\supset
\cw_{N+1}=0$$ such that the quotients $\cw_k/\cw_{k+1}$ have a
structure of $\co_X$-modules.

Then $\cv$ is soft.
\end{proposition}
{\bf Proof.} Let $Z$ be any closed subset of $X$. We have to show
that any section of $\cv$ over $Z$ extends to a section over $X$.
By Lemma \ref{sh-12}(1) it is enough to check that for any open
imbedding $j\colon U\inj X$ one has  $H^i(X,j_!j^*\cv)=0$ for
$i>0$. Observe that $j_!$ and $j^*$ are exact functors (see e.g.
2.5.4 and 2.3.2 of \cite{kashiwara-schapira} respectively). So we
have a filtration
$$j_!j^*\cv=j_!j^*\cw_0\supset j_!j^*\cw_1 \supset \dots \supset
j_!j^*\cw_N\supset j_!j^*\cw_{N+1}=0.$$ By induction and the long
exact sequence, it is enough to check that
$H^i(X,j_!j^*\cw_k/j_!j^*\cw_{k+1})=0$ for $i>0$. But since the
functor $j_!j^*$ is exact we have
$j_!j^*\cw_k/j_!j^*\cw_{k+1}=j_!j^*(\cw_k/\cw_{k+1})$. Now the
result follows from Lemma \ref{sh-12}(3). \qed

\subsection{Valuation theory.}\label{valuations} In this subsection we remind some
facts from the theory of valuations. Let $V$ be a finite
dimensional real vector space, $n=\dim V$. Let $\ck(V)$ denote the
class of all convex compact subsets of $V$. Equipped with the
Hausdorff metric, the space $\ck(V)$ is a locally compact space.

\begin{definition}\label{val-1}
a) A function $\phi :{\cal K}(V) \to \CC$ is called a valuation if
for any $K_1, \, K_2 \in {\cal K}(V)$ such that their union is
also convex one has
$$\phi(K_1 \cup K_2)= \phi(K_1) +\phi(K_2) -\phi(K_1 \cap K_2).$$

b) A valuation $\phi$ is called continuous if it is continuous
with respect to the Hausdorff metric on ${\cal K}(V)$.
\end{definition}

Let us denote by $Val(V)$ the space of translation invariant
continuous valuations on $\ck(V)$. Equipped with the topology of
uniform convergence on compact subsets of $\ck(V)$ the space
$Val(V)$ becomes a Banach space (see e.g. Lemma A.4 in
\cite{alesker-poly}).

\begin{definition}\label{val-2} Let $k$ be a real number.
A valuation $\phi$ is called $k$-homogeneous if for every convex
compact set $K$ and for every scalar $\lam >0$
$$\phi(\lam K)=\lam ^k \phi (K).$$
\end{definition}
Let us denote by $Val _k(V)$ the space of $k$-homogeneous
translation invariant continuous valuations.
\begin{theorem}[McMullen \cite{mcmullen-euler}]\label{val-3}
$$Val(V)=\bigoplus_{k=0}^{n} Val_k(V),$$
where $n=\dim V$.
\end{theorem}
Note in particular that the degree of homogeneity is an integer
between 0 and $n=\dim V$. It is known that $Val_0(V)$ is
one-dimensional and is spanned by the Euler characteristic $\chi$,
and $Val_n(V)$ is also one-dimensional and it is spanned by a
Lebesgue measure \cite{hadwiger-book}. The space $Val_n(V)$ is
also denoted by $|\wedge V^*|$ (the space of complex valued
Lebesgue measures on $V$). One has further decomposition with
respect to parity:
$$Val_k(V)=Val_k^{ev}(V)\oplus Val_k^{odd}(V),$$
where $Val_k^{ev}(V)$ is the subspace of even valuations ($\phi$
is called even if $\phi(-K)=\phi(K)$ for every $K\in {\cal
K}(V)$), and $Val_k^{odd}(V)$ is the subspace of odd valuations
($\phi$ is called odd if $\phi (-K)=-\phi(K)$ for every $K\in
{\cal K}(V)$). The Irreducibility Theorem is as follows.
\begin{theorem}[Irreducibility Theorem
\cite{alesker-gafa}]\label{irr} The natural representation of the
group $GL(V)$ on each space $Val_k^{ev}(V)$ and $Val_k^{odd}(V)$
is irreducible.
\end{theorem}
In this theorem, by the natural representation one means the
action of $g\in GL(V)$ on $\phi\in Val(V)$ as
$(g\phi)(K)=\phi(g^{-1}K)$ for every $K\in \ck(V)$. The subspace
of smooth valuations with respect to this action in sense of
Definition \ref{rep-1} is denoted by $Val^{sm}(V)$.

In \cite{alesker-I} we have introduced the notion of a smooth
valuation on a linear space $V$. Let us remind this notion. Let us
denote by $CV(V)$ the space of continuous valuations on $V$.
Equipped with the topology of uniform convergence on compact
subsets of $\ck(V)$, $CV(V)$ becomes a Fr\'echet space. Let $\qvv$
denote the space of continuous valuations on $V$ which satisfy the
following additional property:
the map given by $K\mapsto \phi( tK +x)$ is a continuous map $\ck
(V) \to C^{n}([0,1]\times V)$. We call such valuations {\itshape
quasi-smooth}.

In the space $\qvv$ we have the natural linear topology defined as
follows. Fix a compact subset $G\subset V$. Define a seminorm on
$\qvv$
$$||\phi||_G:=sup \{
||\phi(tK+x)||_{C^{n}([0,1]\times G)}|\, K\in \ck(V),\, K\subset
G\}.$$ Note that the seminorm $||\cdot||_G$ is finite. One easily
checks the following claim.
\begin{claim}\label{val-5}
Equipped with the topology defined by this family of seminorms the
space $\qvv$ is a Fr\'echet space.
\end{claim}

Note also that the natural representation of the group $Aff(V)$ of
affine transformations of $V$ in the space $\qvv$ is continuous.
We will denote by $\svv$ the subspace of $Aff(V)$-smooth vectors
in $\qvv$. It is a Fr\'echet space.

\begin{definition}\label{val-6}
Elements of $SV(V)$ are called {\itshape smooth} valuations on the
linear space $V$.
\end{definition}

Let us remind notions of characteristic and normal cycle of a
convex compact set $K\in \ck(V)$. Clearly $T^*V=V\times V^*$. Let
$K\in \ck(V)$. Let $x\in K$.
\begin{definition}\label{val-7}
A tangent cone to $K$ at $x$ is a set denoted by $T_xK$ which is
equal to the closure of the set $\{y\in V|\exists \eps>0\, \,
x+\eps y\in K\}$.
\end{definition}
It is easy to see that $T_xK$ is a closed convex cone.
\begin{definition}\label{val-8}
A normal cone to $K$ at $x$ is the set
$$Nor_xK:=\{y\in V^*| \,\, y(x)\geq 0 \forall x\in T_xK\}.$$
\end{definition}
Thus $Nor_xK$ is also a closed convex cone.
\begin{definition}\label{val-9}
Let $K\in \ck(V)$. The {\itshape characteristic cycle} of K is the
set $$CC(K):=\cup_{x\in K}Nor_x(K).$$
\end{definition}
\begin{remark}\label{remark-cc-nc}
The notion of the characteristic cycle is not new. First an almost
equivalent notion of normal cycle (see below) was introduced by
Wintgen \cite{wintgen}, and then studied further by Z\"ahle
\cite{zahle87} by the tools of geometric measure theory.
Characteristic cycles of subanalytic sets of real analytic
manifolds were introduced by Kashiwara (see
\cite{kashiwara-schapira}, Chapter 9) using the tools of the sheaf
theory, and independently by J. Fu \cite{fu-94} using rather
different tools of geometric measure theory. Below we will remind
an elementary definition of characteristic cycle of a
differentiable subpolyhedron of a smooth manifold. This elementary
approach will be sufficient for the purposes of this article.
\end{remark}

For a linear space $W$ let us denote by $\PP_+(W)$ the manifold of
oriented lines in $W$ passing through the origin. Similarly for a
vector bundle $E$ over a manifold $X$ let us denote by $\PP_+(E)$
the vector bundle over $X$ whose fiber over a point $x\in X$ is
equal to $\PP_+(E_x)$ where $E_x$ is the fiber of $E$ over $X$.

It is easy to see that $CC(K)$ is a closed $n$-dimensional subset
of $T^*V=V\times V^*$ invariant with respect to the multiplication
by non-negative numbers acting on the second factor.
\def\ucc{\underline{CC}}
\def\tcc{\tilde{CC}}
Sometimes we will also use the following notation. Let
$\underline{0}$ denote the zero section of $T^*V$, i.e.
$\underline{0}=V\times\{0\}$. Set
\begin{eqnarray*}
\ucc(K):=CC(K)\backslash\underline{0},\\
\tcc(K):=\ucc/\RR_{>0}.
\end{eqnarray*}
Thus $\tcc(K)\subset \PP_+(T^*V)$. Let us denote by $N(K)$ the
image of $\tcc(K)$ under the involution on $\PP_+(T^*V)$ of the
change of an orientation of a line. $N(K)$ is called the {\itshape
normal cycle} of $K$.

In this article for a manifold $Y$ we denote by $\Ome^k:=\wedge^k
T^*Y$ the vector bundle of $k$-forms over $Y$. Usually it will be
clear from the context which manifold is meant.

Let us denote by $$p:T^*V\to V$$ the canonical projection. Let us
denote by $o$ the orientation bundle of $V$. Note that a choice of
orientation on $V$ induces canonically an orientation on $CC(K)$
and $N(K)$ for any $K\in \ck(V)$. Let us denote by $\tilde
C^1(T^*V,\Ome^n\otimes p^*o)$ the space of $C^1$-smooth sections
of $\Ome^n\otimes p^*o$ over $T^*V$ such that the restriction of
$p$ to the support of this section is proper.
\begin{theorem}[\cite{part3}]\label{val-10}
For any $\ome\in\tilde C^1(T^*V,\Ome^n\otimes p^*o)$ the map
$\ck(V)\to \CC$ given by $K\mapsto \int_{CC(K)}\ome$ defines a
continuous valuation on $\ck(V)$.
\end{theorem}

\begin{corollary}\label{val-11}
For any $\eta\in C^1(\PP_+(T^*V),\Ome^{n-1}\otimes p^*o)$  the map
$\ck(V)\to \CC$ given by $K\mapsto \int_{N(K)}\eta$ defines a
continuous valuation on $\ck(V)$.
\end{corollary}
We will also need the following statement.
\begin{theorem}[\cite{part3}]\label{val-12}
The map $\ck(V)\times \left(C^1(V,|\ome_V|)\oplus
C^1(\PP_+(V^*),\Ome^{n-1}\otimes p^*o)\right)\to \CC$ given by
$$(K,(\ome,\eta))\mapsto \int_K\ome +\int_{N(K)}\eta$$
is continuous.
\end{theorem}
\begin{corollary}[\cite{alesker-I}, Corollary 5.1.7]\label{val-13}
(i) The map $C^1(V,|\ome_V|)\oplus
C^1(\PP_+(V^*),\Ome^{n-1}\otimes p^*o)\to CV(V)$ given by
$(\ome,\eta)\mapsto[K\mapsto\int_K\ome +\int_{N(K)}\eta]$ is
continuous.

(ii) For any compact set $G\subset V$ the exists a larger compact
set $\tilde G\subset V$ and a constant $C=C(G)$ such that for any
$(\ome,\eta)\in C^1(V,|\ome_V|)\oplus
C^1(\PP_+(V^*),\Ome^{n-1}\otimes p^*o)$ one has
$$\sup_{K\subset G,K\in\ck(V)}|\int_K\ome+\int_{N(K)}\eta|\leq
C(||\ome||_{C^1(\tilde G)}+||\eta||_{C^1(p^{-1}\tilde G)}).$$
\end{corollary}
\begin{proposition}[\cite{alesker-I}, Proposition 5.1.8]\label{val-14}

(i) For any $$(\ome,\eta)\in C^\infty(V,|\ome_V|)\oplus
C^\infty(\PP_+(V^*),\Ome^{n-1}\otimes p^*o)$$ the valuation
$[K\mapsto\int_K\ome+ \int_{N(K)}\eta]$ is smooth, i.e. belongs to
$SV(V)$.

(ii) The induced map $$C^\infty(V,|\ome_V|)\oplus
C^\infty(\PP_+(V^*),\Ome^{n-1}\otimes p^*o)\to SV(V)$$ is
continuous.
\end{proposition}
\begin{theorem}[\cite{alesker-I}, Theorem 5.2.2]\label{val-15}
The map
$$C^\infty(V,|\ome_V|)\oplus
C^\infty(\PP_+({T^*V}),\Ome^{n-1}\otimes p^*o)\to SV(V)$$ is onto.
\end{theorem}

In \cite{alesker-I} we have defined a decreasing filtration
$W_\bullet$ by closed subspaces on $SV(V)$:
$$SV(V)=W_0\supset W_1\supset \dots\supset W_n.$$
Here
 $$W_i:= \{\phi \in \svv|\, \frac{d^k}{dt^k}
  \phi(tK+x)\big |_{t=0}=0 \,
 \forall k<i,\, \forall K\in \ck(V),\, \forall x\in V\}.$$
 It is clear that $W_i$ are $Aff(V)$-invariant closed subspaces of
 $\svv$.
Obviously $\svv =W_0\supset W_1\supset \dots .$
\begin{proposition}[\cite{alesker-I}, Proposition 3.1.1.]\label{val-16}
$$W_{n+1}=0.$$
\end{proposition}
\begin{proposition}[\cite{alesker-I}, Proposition 3.1.2]\label{val-17}
$W_n$ coincides with the space of smooth densities on $V$.
\end{proposition}

Let us denote by $Val(TV)$ the (infinite dimensional) vector
bundle over $V$ whose fiber over $x\in V$ is equal to the space of
translation invariant $GL(V)$-smooth valuations on the tangent
space $T_xV=V$. Similarly we can define the vector bundle
$Val_k(TV)$ of $k$-homogeneous smooth translation invariant
valuations. Clearly
$C^\infty(V,Val_k(TV))=C^\infty(V,Val_k^{sm}(V))$ where the last
space denotes the space of infinitely smooth functions on $V$ with
values in the Fr\'echet space $Val^{sm}_k(V)$.
\begin{theorem}[\cite{alesker-I}, Proposition 3.1.5.]\label{val-18} There exists
a canonical isomorphism of Fr\'echet spaces of the associated
graded space $gr_W SV(V):=\bigoplus_{i=0}^nW_i/W_{i+1}$ and
$C^\infty(V,Val(TV))$.
\end{theorem}
Remind also the construction of this isomorphism. More precisely
there is an isomorphism
$$I_i:W_i/W_{i+1}\tilde \to C^\infty(V,Val_i^{sm}(V)).$$
The map $I_i$ is defined as follows. For $\phi\in W_i,x\in V, K\in
\ck(V) $
\begin{eqnarray}\label{Ii}
(I_i\phi)(x,K):=\lim_{r\to +0}\frac{\phi(rK+x)}{r^i}.
\end{eqnarray}

Now let us describe the filtration $W_\bullet$ in terms of
integration with respect to the characteristic cycle following
\cite{alesker-I}. Let us start with some general remarks.

Let $X$ be a smooth manifold. Let $p:P\to X$ be a smooth bundle.
Let $\Omega^N(P)$ be the vector bundle over $P$ of $N$-forms. Let
us introduce a filtration of $\Omega^N(P)$ by vector subbundles
$W_i(P)$ as follows. For every $y\in P$ set
\begin{multline*} (W_i(P))_y:=\{\omega\in \wedge^NT_y^*P \big|\,
\omega|_F\equiv 0 \mbox{ for all } F\in Gr_N(T_yP)\\
\mbox{ with } \dim(F\cap T_y(p^{-1}p(y)))>N-i\}.\end{multline*}

Clearly we have
$$\Omega^N(P)=W_0(P)\supset W_1(P)\supset \dots \supset W_N(P)\supset
W_{N+1}(P)=0.$$ Let us discuss this filtration in greater detail
following \cite{alesker-I}.

Let us make some elementary observations from linear algebra.
\def\wel{W(L,E)}
Let $L$ be a finite dimensional vector space. Let $E\subset L$ be
a linear subspace. For a non-negative integer $i$ set
$$\wel_i:=\{\omega\in \wedge^NL^*\big|\, \omega|_F\equiv 0\mbox{ for
all } F\subset L \mbox{ with } \dim(F\cap E)>N-i\}.$$ Clearly
$$\wedge^NL^*=\wel_0\supset \wel_1\supset \dots \supset \wel_N\supset
\wel_{N+1}=0.$$
\begin{lemma}[\cite{alesker-I}, Lemma 5.2.3] \label{val-19}
There exists canonical isomorphism of vector spaces
$$\wel_i/\wel_{i+1}=\wedge^{N-i}E^*\otimes \wedge^i(L/E)^*.$$
\end{lemma}

Let us apply this construction in the context of integration with
respect to the characteristic cycle. Let $X$ be a smooth manifold
of dimension $n$. Let $P:=T^*X$ be the cotangent bundle. Let
$p:P\to X$ be the canonical projection. Let us denote by $o$ the
orientation bundle on $X$.  The above construction gives a
filtration of $\Omega^n(P)$ by subbundles
$$\Omega^n(P)=W_0(\Omega^n(P))\supset\dots\supset W_n(\Omega^n(P)).$$
Twisting this filtration by $p^*o$ we get a filtration of
$\Omega^n(P)\otimes p^*o$ by subbundles denoted by
$W_i(\Omega^n(P)\otimes p^*o)$.

Let us denote by $\tilde C^\infty(P,W_i(\Omega^n\otimes p^*o))$
the space of infinitely smooth sections of the bundle
$W_i(\Omega^n\otimes p^*o)$ such that the restriction of the
projection $p$ to the support of these sections is proper. The
next result is a trivial reformulation of Proposition 5.2.5 from
\cite{alesker-I}.
\begin{theorem}\label{val-20}
Consider the map $\Xi: \tilde C^\infty(P,\Omega^n\otimes p^*o)\to
SV(V)$ given by $$(\Xi(\omega))(K)=\int_{CC(K)}\omega.$$ This map
is surjective, and moreover for every $i=0,1,\dots,n$ the map
 $\Xi$ maps $\tilde C^\infty(P,W_i(\Omega^n)\otimes p^*o)$ onto
$W_i$ surjectively.
\end{theorem}

\section{Differentiable polyhedra, finitely additive
measures, and smooth valuations.}\label{pomeva} In  Subsection
\ref{poly} we discuss the notion of the differentiable polyhedron.
In Subsection \ref{fa} we discuss a combinatorial notion of
finitely additive measure on a family of sets which is not
necessarily closed under finite intersections and unions but
satisfies some other assumptions. In Subsection \ref{sheaf1} we
introduce a notion of finitely additive measure on the class of
simple differentiable subpolyhedra of a smooth manifold. Finally
in Subsection \ref{sheaf-sm} we introduce the main object of this
article, namely the notion of a smooth valuation on a manifold.
\subsection{Differentiable polyhedra.} \label{poly}
 We remind the definition and basic properties of
differentiable polyhedra. The exposition in the beginning this
subsection (up to Lemma \ref{poly-6})  follows very closely
\cite{allendoerfer-weil}.

A {\itshape convex angle} in $\RR^n$ is a set defined by finitely
many inequalities $\{x|\, <\xi_\nu,x>\geq 0 , \, 0\leq \nu\leq
N\}$. Note that a convex angle is a convex cone in particular. We
say that a convex angle $C$ is of type $r$ if it contains an
$r$-dimensional linear subspace and does not contain linear
subspaces of larger dimensions.


Let $P^n$ be a compact connected topological space for which there
has been given a covering by open subsets $\Ome_i$ and a
homeomorphic mapping $\phi_i$  of each $\Ome_i$ onto an
$n$-dimensional convex angle $C_i$ (which may be $\RR^n$). $P^n$
is called an $n$-dimensional {\itshape differentiable polyhedron}
if the maps $\phi_i\phi_k^{-1}$ are of class $C^\infty$ on the
domain of their definition.

A {\itshape differentiable cell} is, by definition, a
differentiable polyhedron  which is diffeomorphic of class
$C^\infty$ with a convex compact polyhedron in $\RR^n$.

Let $P^n$ be a differentiable polyhedron. For any point $z\in P^n$
one defines the tangent space $\tilde T_zP^n$ to $P^n$ at $z$ in
the obvious way. $\tilde T_zP^n$ is a linear space. The (tangent)
angle of $P^n$ at $z$ is the subset of $\tilde T_zP^n$ consisting
of those $v\in \tilde T_zP^n$ such that there exists a
$C^\infty$-smooth map $\gamma:[0,1]\to P^n$ such that
$\gamma(0)=z,\, \gamma'(0)=v$. We will denote it by $T_zP^n$. It
is clear that $T_zP^n$ is a convex angle in $\tilde T_zP^n$.


If $C$ is the tangent angle of $P^n$ at a point $z$ then $z$ has a
neighborhood homeomorphic to $C$. If $C$ is of type $r$ then we
say that $z$ is of type $r$ in $P^n$. Points of type $n$ in $P^n$
are called {\itshape inner} points of $P^n$. Points of type at
most $r$ (where $0\leq r\leq n$) form  a compact subset of $P^n$.
The set of inner points of $P^n$ will be called (relative)
interior of $P^n$ and will be denoted by $int P^n$.

\begin{definition}\label{poly-1}
 A {\itshape regular differentiable subpolyhedron} $Q^p$ in $P^n$, is the
one-to-one image of a differentiable polyhedron $Q^p_0$ in $P^n$
provided that this map is of class $C^\infty$ and its differential
is injective at every point.
\end{definition}

\begin{definition}\label{poly-2}
A finite set of distinct regular subpolyhedra $Q^r_\rho$ of $P^n$
form a {\itshape subdivision} $\cd$ of $P^n$ if the following
conditions are satisfied:

(1) each point of $P^n$ is an inner point of at least one
$Q^r_\rho$ in $\cd$;

(2) if $Q^r_\rho$ and $Q^s_\sigma$ in $\cd$ are such that there is
an inner point of $Q^r_\rho$ contained in $Q^s_\sigma$, then
$Q^r_\rho\subset Q^s_\sigma$.
\end{definition}

>From condition (2) it follows that no two differentiable polyhedra
in $\cd$ can have an inner point in common unless they coincide.

The following result was proved in \cite{allendoerfer-weil}, Lemma
7.
\begin{lemma}\label{poly-3}
If $Q^r$ is a differentiable polyhedron in a subdivision $\cd$ of
$P^n$, all inner points of $Q^r$ have the same type in $P^n$.
\end{lemma}

\begin{lemma}\label{poly-4}
Let $\cd$ be a subdivision of $P^n$. Let $A$ and $B$ be two
subsets of $P^n$ which are unions of finitely many elements of the
subdivision $\cd$. Then $A\cap B$ is also a union of finitely many
elements of $\cd$.
\end{lemma}
{\bf Proof.} It is enough to prove the lemma under the assumption
that $A$ and $B$ are elements of $\cd$. Assume that $z\in A\cap
B$. Then there is a unique cell  $P^r_\lam\in \cd$ such that $z$
belongs to its interior. Then by part (2) of Definition
\ref{poly-2} $P^r_\lam\subset A$ and similarly $P_\lam^r\subset
B$. Hence $P^r_\lam \subset A\cap B$. The result follows. \qed

\begin{definition}\label{poly-5}
 A subdivision $\cd '$ of $P^n$ is called
 a {\itshape refinement}  of
a subdivision $\cd$ of $P^n$ if for any differentiable polyhedron
$P^r_\lam$ in $\cd$ all differentiable polyhedra of $\cd'$
contained in $P^r_\lam$ form a subdivision of $P^r_\lam$.
\end{definition}

Lemma \ref{poly-3} implies that if a differentiable polyhedron
$Q^r$, in a subdivision $\cd$ of $P^n$, contains at least one
inner point of $P^n$ then all inner points of $Q^r$ are inner
points of $P^n$; $Q^r$ is called an {\itshape inner polyhedron} of
the subdivision $\cd$.

\begin{lemma}[\cite{allendoerfer-weil}, Lemma 8]\label{poly-6}
Let $\cd$ be a subdivision of $P^n$ and let $z$ be any point of
$P^n$. Then the tangent angles at $z$ of those differentiable
polyhedra in $\cd$ which contain $z$ form a subdivision of the
tangent angle of $P^n$ at $z$. Moreover the inner angles (i.e.
those of maximal dimension) in the latter subdivision are the
angles of the inner polyhedra in $\cd$ which contain $z$.
\end{lemma}

A {\itshape differentiable cell } is a differentiable polyhedron
diffeomorphic to a convex compact polytope. We now define a
{\itshape cellular} subdivision of a differentiable polyhedron
$P^n$ as a subdivision $\cd$, every polyhedron $Z^r_\rho$ in which
is a differentiable cell.

\begin{definition}\label{poly-7}
(1) A differentiable polyhedron $P^n$ is called {\itshape simple}
if every point $z\in P^n$ has a neighborhood diffeomorphic to
$\RR^r \times \RR^{n-r}_{\geq 0}$ for some $0\leq r\leq n$.

(2) A subdivision $\cd$ of a differentiable polyhedron $P^n$ is
called {\itshape simple} if any element of $\cd$ is simple.

(3) A {\itshape triangulation} of a differentiable polyhedron
$P^n$ is a subdivision every element of which is diffeomorphic to
a simplex.
\end{definition}

The following result is well known.
\begin{proposition}\label{poly-8}
Every simple polyhedron admits a triangulation.
\end{proposition}


\begin{definition}\label{poly-9}
Let $\cd=\{P_\lam\}$ and $\cd'=\{P_\nu'\}$ be two subdivisions of
a differentiable polyhedron $P$. We say that $\cd$ and $\cd'$ are
transversal to each other if for any $P_\lam\in\cd$, any
$P_\nu'\in \cd'$, and any $x\in P_\lam\cap P_\nu'$ the maximal
linear subspaces contained in the cones $T_xP_\lam$ and
$T_xP_\nu'$ intersect transversally.
\end{definition}

\begin{lemma}\label{poly-10}
Let $X^{(1)}$ and $X^{(2)}$ be two regular differentiable
subpolyhedra of a smooth $n$-manifold $M$. Assume that for any
$x\in X^{(1)}\cap X^{(2)}$ the maximal linear subspaces contained
in the cones $T_xX^{(1)}$ and $T_xX^{(2)}$ intersect
transversally. Then $X^{(1)}\cap X^{(2)}$ is a differentiable
polyhedron. Moreover if $X^{(1)}$ and $X^{(2)}$ are simple then
$X^{(1)}\cap X^{(2)}$ is also simple.
\end{lemma}
{\bf Proof.} Fix $x\in X^{(1)}\cap X^{(2)}$.  Let $x$ has type
$p_1$ in $X^{(1)}$ and type $p_2$ in $X^{(2)}$. Then there exist
$C^\infty$-smooth real  valued functions
$f_1^{(k)},\dots,f_{N_k}^{(k)},\, k=1,2$, such that

(1) for $k=1,2$ for each $j>n-p_k$ the function $f_j^{(k)}$ is a
linear combination with constant coefficients of $f_l^{(k)}$'s
with $l\leq n-p_k$;

(2) in a neighborhood of $x$
$$X^{(k)}=\{f_j^{(k)}\geq 0|\, 1\leq j\leq N_k\};$$

(3)
$df_1^{(1)}|_x,\dots,df_{n-p_1}^{(1)}|_x;df_1^{(2)}|_x,\dots,df_{n-p_2}^{(2)}|_x$
are linearly independent.

Let $q:=n-((n-p_1)+(n-p_2)).$ Let us choose $C^\infty$-smooth
functions $g_1,\dots ,g_q$ such that
$$df_1^{(1)}|_x,\dots,df_{n-p_1}^{(1)}|_x;df_1^{(2)}|_x,\dots,df_{n-p_2}^{(2)}|_x;
dg_1|_x,\dots,dg_q|_x$$ form a basis of $T_x^*M$. Then the
sequence of functions
$$f_1^{(1)},\dots,f_{n-p_1}^{(1)};f_1^{(2)},\dots,f_{n-p_2}^{(2)};g_1,\dots,g_q$$
form a coordinate system in a neighborhood of $x$. It is clear
that in this coordinate system $X^{(1)},X^{(2)},X^{(1)}\cap
X^{(2)}$ are defined by finite systems of linear inequalities, and
hence they are convex angles. The last part of the lemma also
follows. \qed

\begin{proposition}\label{poly-11}
Let $\cd=\{P_\lam\}$ and $\cd'=\{P_\nu'\}$ be two transversal
subdivisions of a polyhedron $P$. Let $\cd \cap \cd':=\{P_\lam\cap
P_\nu'|\, P_\lam\in \cd,\, P_\nu'\in \cd'\}$. Then $\cd\cap\cd'$
is also a subdivision of $P$. Moreover it refines both $\cd$ and
$\cd'$.
\end{proposition}
 To prove this proposition we need first of all the
following lemma.

\begin{lemma}\label{poly-12}
Let $M$ be a smooth manifold. Let $P$ and $Q$ be two regular
differentiable subpolyhedra of $M$. Assume that for any $x\in
P\cap Q$ the maximal linear subspaces contained in the cones
$T_xP$ and $T_x Q$ intersect transversally. Then $int (P\cap
Q)=int P\cap int Q$.
\end{lemma}
{\bf Proof.} By Lemma \ref{poly-10} $P\cap Q$ is a differentiable
polyhedron. Fix $z\in P\cap Q$. Let $z$ has type $p$ in $P$ and
type $q$ in $Q$. Consider the tangent space $T_zM$ to $M$ at $z$
and tangent angles $P_1$ and $Q_1$ to $P$ and $Q$ respectively at
$z$. Then $P_1,\, Q_1\subset T_zM$ are convex angles. $P_1$
contains a $p$-dimensional linear subspace $A$ such that the image
of $P_1$ in $T_zM/A$ is a {\itshape cornered} convex angle (i.e.
it does not contain any non-zero linear subspace). Similarly $Q_1$
contains a $q$-dimensional linear subspace $B$ such that the image
of $Q_1$ in $T_zM/B$ is a cornered convex angle. It follows from
the assumptions of the lemma that $A$ and $B$ are transversal to
each other. Put $C:=A\cap B$. Choose $A'$ a complement of $C$ in
$A$, and $B'$ a complement of $C$ in $B$. Then there exist
cornered convex angles $R\subset B'$ and $S\subset A'$ such that
$$P_1=A\times R=C\times A'\times R,\, Q_1=B\times S=C\times
S\times B'.$$ Then $P_1\cap Q_1=C\times S\times R$. This is also a
convex angle and $int(P_1\cap Q_1)=C\times int S\times int R=int
P_1\cap int Q_1. $ It is easy to see that in a small neighborhood
of $z$ the polyhedron $P\cap Q$ is diffeomorphic to $P_1\cap Q_1$.
This implies the lemma.
 \qed

{\bf Proof} of Proposition \ref{poly-11}.
 Now let us check that $\cd\cap \cd'$ is indeed a
subdivision of $P$. Fix any $z\in P$. Then by Lemma 2.1.12 there
exits $P_\lam\in \cd$ and $P_\nu'\in \cd'$ such that $z\in int
P_\lam$ and $z\in int P_\nu'$. Hence $z\in int P_\lam\cap int
P_\nu'=int (P_\lam\cap P_\nu')$.

It remains to check condition (2) of Definition \ref{poly-2}.
Assume that $z\in int (P_\lam \cap P_\nu')=int P_\lam\cap int
P_\nu'$ and $z\in P_s\cap P_t'$. Then it follows that $P_\lam
\subset P_s$ and $P_\nu'\subset P_t'$. Hence $P_\lam \cap
P_\nu'\subset P_s\cap P_t'$. \qed

\begin{definition}\label{poly-13}
Let $\cd$ be a  subdivision of a differentiable polyhedron $P$.
Let $\{U_\alp\}$ be an open covering of $P$. We say that $\cd$ is
subordinate to $\{U_\alp\}$ if any element of $\cd$ is contained
in at least one element of the covering $\{U_\alp\}$.
\end{definition}
\begin{lemma}\label{poly-14}
Let $\cd$ be a cellular subdivision of a differentiable polyhedron
$P$. Let $\{U_\alp\}$ be an open covering of $P$. Then there
exists a refinement $\cd'$ of $\cd$ which is a triangulation of
$P$ and subordinate to $\{U_\alp\}$.
\end{lemma}
{\bf Proof.} Assume that we have constructed a triangulation of
each element of $\cd$ of dimension less than $r$ subordinate to
$\{U_\alp\}$. Let us fix a cell $Q\in \cd$ of dimension $r$ and
let us construct its subdivision which extends the subdivision of
the boundary of $Q$ obtained previously and which is subordinate
to $\{U_\alp\}$. Let us fix a point $x\in int Q$. For any cell $R$
contained in the boundary of $Q$ and belonging to the subdivision
constructed previously, let us consider the cone over $R$ with
vertex at $x$. All such cones form a subdivision of $Q$.

Now we are reduced to the following situation. Given a convex
compact polytope $S$ of dimension $r$ and given its
$(r-1)$-dimensional face $F\subset S$ which is a simplex. Given an
open covering $\{U_\alp\}$ of $S$ such that $F$ is contained in at
least one of the elements of the covering.  We have to find a
triangulation of $S$ subordinate to $\{U_\alp\}$ such that $F$ is
one of the elements of this subdivision. But this problem in the
affine space can be solved easily. \qed

\begin{lemma}\label{poly-15}
Every differentiable polyhedron can be regularly imbedded into a
smooth compact manifold.
\end{lemma}
{\bf Proof.} We can choose a finite open covering
$\{U_\alp\}_{\alp=1}^N$ of $P$, open sets $\{V_\alp\}_{\alp=1}^N$
such that $\bar{U_\alp}\subset V_\alp$, and there exist
diffeomorphisms $\phi_\alp$ of $V_\alp$ onto a convex angle in
$\RR^n$. Let $l_\alp^1,\dots, l_\alp^n$ be the corresponding
coordinate functions on $V_\alp$. Let us choose a smooth partition
of unity $\{\phi_\alp\}$ such that $\phi_\alp \equiv 1 $ on
$\bar{U_\alp}$, $supp ( \phi_\alp)\subset V_\alp$, and $\sum_\alp
\phi_\alp\equiv 1$. Then the collection of functions $\{\phi_\alp
l_\alp^j\}$ defines an immersion of $P$ into $\RR^{nN}$. Indeed
let us fix $x_0\in P$. There exists $\alp_0$ such that $x_0\in
U_{\alp_0}$. Then for any $x\in U_{\alp_0}$
$$(\phi_{\alp_0}l^j_{\alp_0})(x)=l^j_{\alp_0}(x).$$
hence the functions $\{\phi_{\alp_0}l^j_{\alp_0}\}$ define an
imbedding of $U_{\alp_0}$ to $\RR^N$. Hence all the functions
$\{\phi_\alp l^j_\alp\}$ define an imbedding of $U_{\alp_0}$ to
$\RR^{nN}$. Thus the functions $\{\phi_\alp l^j_\alp\}$ define an
immersion of $P$ to $\RR^{nN}$.

Now let us assume that we have an immersion $f\colon P\to \RR^M$.
Let us construct an imbedding $g\colon P\to \RR^{M'}$. Note that
the fibers of $f$ are discrete sets. Since the space $P$ is
compact the cardinality of fibers of $f$ is uniformly bounded. Fix
a point $y_0\in \RR^M$ lying in the image of $f$. Let
$f^{-1}(y_0)=\{x_1,\dots,x_k\}\subset P$. One can choose a smooth
function $g_{y_0}\colon P\to \RR$ such that
$$g_{y_0}(x_i)\ne g_{y_0}(x_j) \mbox{ for } i\ne j.$$
It is clear that there exists a neighborhood $\co$ of $y_0$ such
that for any $y\in \co$ the function $g_{y_0}$ takes different
values on points from $f^{-1}(y)$. Choosing a finite covering
$\{\co_\beta\}$ of $Im f$ by such subsets we construct smooth
functions $\{g_\beta\}_{\beta=1}^k$ such that for each
$\beta=1,\dots,k$, $g_\beta\colon P\to \RR$, and for any $y\in Im
f$ there exists $\gamma=1,\dots,k$ such that $g_\gamma$ takes
different values on points from $f^{-1}(y)$.

Consider the map $$g:=(f,g_1,\dots,g_k)\colon P\to \RR^M\times
\RR^k=\RR^{M+k}.$$ Obviously this map is an imbedding. Since
$\RR^{M+k}$ can be imbedded as an open subset into the sphere
$S^{M+k}$ the result follows. \qed

\begin{proposition}\label{poly-16}
Let $M$ be a compact smooth manifold. Let $P\subset M$ be a
differentiable polyhedron. Let $\cd$ be a subdivision of $P$. Let
$T$ be a subdivision of $M$. Then the set of
$C^{\infty}$-diffeomorphisms $f$ of $M$ such that for each
$T_\lam\in T$ its image $f(T_\lam)$ is transversal to each $P_\nu
\in \cd$, is open and dense in the group $Diff(M)$ of all
$C^\infty$-diffeomorphisms of $M$.
\end{proposition}
{\bf Proof.} The openness is obvious. Let us prove the density.
Clearly it is enough to prove that in any neighborhood of the
identity diffeomorphism of $M$ there is a transformation we need.
Let $n:=\dim M$. We can choose  a finite open covering
$\{U_i\}_{i=1}^N$ of $M$ such that there exist open subsets
$\{V_i\}_{i=1}^N$ such that $\bar{U}_i\subset V_i$, and there
exist diffeomorphisms  $f_i:V_i\to D^n$ where $D^n$ denotes the
unit ball in $\RR^n$.  Let $l_i^1,\dots,l_i^n$ be the
corresponding coordinate functionals on $V_i$. Let us fix a
partition of unity $\{\phi_i\}_{i=1}^N$ such that $\phi_i \equiv 1
$ on $\bar{U}_i$, $supp ( \phi_i)\subset V_i$, and $\sum_i \phi_i
\equiv 1$. Then for small enough real numbers $a_{ij}$ the map
$x\mapsto x+\sum_{ij}a_{ij}\phi_i l_i^j$ is a globally defined
diffeomorphism of $M$. Let $A\subset \RR^{nN}$ be a small
neighborhood of 0 in the space of parameters $\{a_{ij}\}$. Thus we
get a map $$\Xi:A\times M\to M.$$ Let us fix $T_\lam \in T$ and
$P_\nu\in \cd$. It is clear that if we restrict $\Xi$ to $A\times
T_\lam$ we get a submersion $$\Xi':A\times T_\lam \to M.$$ In
particular $\Xi'$ is transversal to $P_\nu$. Then by Theorem
10.3.3 of \cite{dubrovin-fomenko-novikov} for $a$ from a dense
subset of $A$ the map
$$\Xi_a':=\Xi'(a,\cdot):T_\lam \to M$$
is transversal to $P_\nu$. (Though in
\cite{dubrovin-fomenko-novikov} this is proved under assumption
that $T_\lam,\, P_\nu$ are closed submanifolds, but the same proof
works when $T_\lam$ and $P_\nu$ are differentiable subpolyhedra.)
\qed

\begin{proposition}\label{poly-17}
Let $P$ be a differentiable polyhedron. Let $\{U_\alp\}$ be an
open covering of $P$. Let $\cd$ be a simple subdivision of $P$.
Then there exists a refinement $\cd'$ of $\cd$ which is simple and
subordinate to $\{U_\alp\}$.
\end{proposition}
{\bf Proof.} Using Lemma \ref{poly-15} let us imbed $P$ into a
smooth compact manifold $M$. Let $\tilde U_\alpha$ be an open
subset of $M$ such that $\tilde U_\alpha \cap P=U_\alpha$.
Consider the open covering of $M$ by $\{U_\alp\}\cup \{M\backslash
P\}$. Since any smooth manifold admits a triangulation, by Lemma
\ref{poly-14} we can choose a triangulation $T$ of $M$ subordinate
to this covering. By Proposition \ref{poly-16} we can choose a
generic diffeomorphism of $M$ close to the identity so that the
image of $T$ is transversal to $\cd$. We may assume that it is $T$
itself. Choosing  $\cd ':= \cd \cap T$ and applying Lemma
\ref{poly-10} and Proposition \ref{poly-11} we prove the
proposition. \qed

\begin{proposition}\label{poly-18}
Let $P$ be a differentiable polyhedron. Let $\cd_1$ and $\cd_2$ be
two subdivisions of $P$. Then there exist subdivisions
$\cd_3,\,\cd',\, \cd''$ of $P$ such that

(1) $\cd'$ is a refinement of $\cd_1$ and $\cd_3$;

(2) $\cd''$ is a refinement of $\cd_2$ and $\cd_3$.

Moreover if $\cd_1$ and $\cd_2$ are simple then $\cd_3,\,\cd'$,
and $\cd''$ can also be chosen simple.
\end{proposition}
{\bf Proof.} Using Lemma \ref{poly-15} let us fix an imbedding of
$P$ into a compact smooth manifold $M$. Fix any triangulation $T$
of $M$. Let $\cd_3$ be the image of $T$ under a generic
diffeomorphism of $M$ (we use Proposition \ref{poly-16}). Then
$\cd_3$ is transversal to $\cd_1$ and $\cd_2$ . Now let us define
$\cd':=\cd_1\cap \cd_3$, $\cd'':=\cd_2\cap \cd_3$. The result now
follows from Lemma \ref{poly-10} and Proposition \ref{poly-11}.
\qed

\begin{definition}\label{poly-19}
Let $\cd$ be a subdivision of $P^n$. Assume that a subset
$X\subset P^n$ admits a presentation as a union $X=\cup_{j=1}^s
P_{\lam_j}$ where $P_{\lam_j}\in \cd$. We say that this
presentation of $X$ is {\itshape reduced} if no one of the
polyhedra in this union is contained in another, i.e.
$P_{\lam_i}\not\subset P_{\lam_j}$ for $i\ne j$.
\end{definition}
\begin{lemma}\label{poly-20}
Let us assume that a subset $X\subset P^n$ has two reduced
decompositions
$$X=A\cup (\cup_jP_{\lam_j})=A\cup (\cup _l Q_{\nu_l})$$
where $A\in \cd$ is a polytope of type $r$ and $P_{\lam_j},\,
Q_{\nu_l}$ are polytopes of type at most $r$. Then
$$\cup_jP_{\lam_j}=\cup _l Q_{\nu_l}.$$
\end{lemma}
{\bf Proof.} Set $B:=\cup_jP_{\lam_j}, \, C:=\cup _l Q_{\nu_l}$.
By symmetry it is enough to prove that $B\subset C$. Let $z\in B$,
say $z\in P_{\lam_1}$. If $z\not \in A$ then $z\in C$. Let us
assume now that $z\in A$. By assumption $P_{\lam_1}\not \subset
A$. Fix any point $w$ from the interior of $P_{\lam_1}$. Then
$w\not\in A$. Hence $w\in C$. Hence $P_{\lam_1} \subset C$. In
particular $z\in C$. \qed

\begin{corollary}\label{poly-21}
Let $\cd$ be a subdivision of $P^n$. Let $X\subset P^n$ be a
subset presentable as a union of some elements of $\cd$. Then $X$
admits a reduced decomposition, and it is unique.
\end{corollary}
{\bf Proof.} The existence of a reduced decomposition is obvious.
Let us prove the uniqueness. Let us denote by $r:=\dim X$. Assume
that we have two reduced decompositions of $X$:
\begin{equation}\label{red}
X=\cup_{j=1}^s A_{\lam_j}=\cup_{l=1}^t B_{\nu_l}.\end{equation}
Take some $B_{\nu_p}$ of dimension $r$. Fix any interior point $z$
of $B_{\nu_p}$. Then $z\in A_{\lam_q}$ for some $\lam_q$. Hence
$B_{\nu_p}\subset A_{\lam_q}$. Since $A_{\lam_q}$ has dimension at
most $r$ we conclude that $B_{\nu_p}=A_{\lam_q}$. By Lemma
\ref{poly-20} we can omit $B_{\nu_p}=A_{\lam_q}$ from the second
equality in (\ref{red}). Continuing this process we prove the
statement. \qed

\subsection{Finitely additive measures.}\label{fa}
In this subsection we will discuss a combinatorial notion of
finitely additive measure on a family of sets which is not
necessarily closed under finite intersections and unions but
satisfies some other assumptions.

\begin{definition}\label{fa-1}
Let $\cs$ be a family of sets which is closed under finite unions
and finite intersections. A functional $\mu\colon \cs \to \CC$ is
called a {\itshape finitely additive measure} if for any $A,B\in
\cs$ one has $$\mu(A\cup B)=\mu(A)+\mu(B)-\mu(A\cap B).$$ It is
easy to see by induction that finitely additive measures satisfy a
stronger inclusion-exclusion property. Namely for any $A_1,\dots ,
A_s \in \cs$ one has
$$\mu(\cup_{i=1}^sA_i)=\sum_{I\subset \{1,\dots, s\},\, I\ne
\emptyset} (-1)^{|I|+1}\mu(\cap_{i\in I} A_i).$$
\end{definition}

Let $\cd=\{A_\lam\}_{\lam \in \Lam}$ be a finite family of subsets
of some set. Assume that we are given a decomposition of the set
of indices $\Lam$ into a disjoint union $$\Lam =\Lam_0
\coprod\Lam_1\coprod \dots \coprod \Lam_n.$$ For $r=0,1,\dots, n$
let us call sets $A_\lam$ with $\lam \in \Lam_r$ the sets of type
$r$. Set $\Lam_{\leq r}:=\cup_{i=0}^r \Lam_i.$ Let $X$ be a finite
union of some of elements of $\cd$. Let us call a presentation
$X=\cup_{j=1}^s A_{\lam_j},\, \lam_j\in \Lam$ {\itshape reduced}
if no set $A_{\lam_j}$ in this presentation is contained in the
other one.

Let us make the following assumptions on $\cd$:

(1) for any sets $A_{\lam_1}$ and $A_{\lam_2}$ from $\cd$ of types
$r_1$ and $r_2$ respectively, their intersection $A_{\lam_1}\cap
A_{\lam_2}$ is a finite union of sets from $\cd$ of types at most
$\min \{r_1,r_2\}$;

(2) if for $\lam_1\ne \lam_2$ the sets $A_{\lam_1}$ and
$A_{\lam_2}$ are of the same type $r$ then $A_{\lam_1}\cap
A_{\lam_2}$ is a finite union of sets from $\cd$ of types strictly
less than $r$;

(3) For every set $X$ as above, a reduced decomposition is unique.

Let us denote by $\ct$ the family of all finite unions of subsets
from  $\cd$. Then clearly under the above assumptions $\ct$ is
closed under finite unions and finite intersections.

Assume we are given a function $m:\Lam \to \CC$. Then we have
\begin{lemma}\label{fa-2}
Under the above assumptions there exists unique finitely additive
measure $\mu$ on $\ct$ (in the sense of Definition \ref{fa-1})
such that for any $\lam \in\Lam$
$$\mu(A_\lam)=m(\lam).$$
\end{lemma}
{\bf Proof.} Let us denote by $\ct_r$ the family of all finite
unions of subsets from $\cd$ of types at most $r$. Then clearly
$\ct_r$ is closed under finite unions and finite intersections.
Moreover we have:
$$\ct_0\subset \ct_1\subset \dots \subset \ct_n=\ct.$$
The construction of the measure $\mu$ on $\ct_r$ will be by
induction on $r$. First let $r=0$. For any $\lam_1, \lam_2\in
\Lam_0,\, \lam_1\ne \lam_2$ we have $A_{\lam_1}\cap
A_{\lam_2}=\emptyset$. Hence any set $X$ from $\ct_0$ has the form
$X= A_{\lam_1}\coprod \dots \coprod A_{\lam_s}$ where all
$\lam_j\in \Lam _0$ and are distinct. Then there is only one way
to define $\mu$ on $\ct_0$, namely $\mu(X)=\sum_{j=1}^s
m(A_{\lam_j})$. Clearly we get a well defined measure $\mu$ on
$\ct_0$.

Assume we have constructed uniquely defined finitely additive
measure $\mu$ on $\ct_{r-1}$. Let us extend it to $\ct_r$ and
prove uniqueness of this extension. Let $X\in \ct_r$. Then $X$ has
a (non-unique) presentation $X=\cup_{j=1}^s A_{\lam_j}$ where
$\lam_j\in \Lam_{\leq r}$ and all $A_{\lam_j}$ are pairwise
distinct. The only way to define $\mu(X)$ is
\begin{equation}\label{addi}
\mu(X):= \sum_{j=1}^s m(\lam_j)+ \sum_{I\subset \{1,\dots,s\},\,
|I|>1} (-1)^{|I|+1} \mu(\cap _{i\in I} A_{\lam_i}).\end{equation}
Note that in this formula the second sum is defined by the
assumption of induction. We have to prove that $\mu$ is well
defined on $\ct_r$ and that it is indeed a finitely additive
measure.

Let us check first that $\mu$ is well defined. Let $X\in \ct_r$.
It is sufficient to show that for any decomposition of a set $X\in
\ct_r$ the expression (\ref{addi}) gives the same value for $\mu$
as for the reduced decomposition (which is unique by the
assumptions on $\cd$). Assume that a decomposition $X=\cup_{j=1}^s
A_{\lam_j}$ is not reduced, say $A_{\lam_1}\supset A_{\lam_2}$.
Then we have:
$$\sum_{I\subset \{1,\dots,s\},I\ne\emptyset}(-1)^{|I|+1} \mu(\cap_{i\in I}
A_{\lam_i})=$$ $$\sum_{I\subset
\{3,\dots,s\},I\ne\emptyset}(-1)^{|I|+1}\mu(\cap_{i\in I}
A_{\lam_i})+\sum_{I\subset
\{3,\dots,s\}}(-1)^{|I|+2}\mu(A_{\lam_1}\cap(\cap_{i\in I}
A_{\lam_i}))+$$ $$ \sum_{I\subset
\{3,\dots,s\}}(-1)^{|I|+2}\left[\mu (A_{\lam_2}\cap(\cap_{i\in I}
A_{\lam_i}))-\mu(A_{\lam_1}\cap A_{\lam_2}\cap(\cap_{i\in I}
A_{\lam_i}))\right].$$ The last sum clearly vanishes. Hence we see
that the set $A_{\lam_2}$ can be omitted. We can continue this
procedure till we get a reduced decomposition of $X$. This proves
that $\mu$ is well defined on $\ct_r$.

It remains to check that $\mu$ is indeed a finitely additive
measure on $\ct_r$. Since $\ct_r$ is closed under finite unions
and finite intersections it is sufficient to check that for any
two sets $X,Y\in \ct_r$ one has $\mu(X\cup
Y)=\mu(X)+\mu(Y)-\mu(X\cap Y).$ Let $X=\cup_{j=1}^s A_{\lam_j},\,
Y=\cup _{l=1}^tB_{\nu_l}$. Let us prove the statement by the
induction in $s$.

Let us assume that $s=1$. Thus $X=A_{\lam_1}$. First let us check
that
\begin{eqnarray}\label{add}
\mu(X\cap Y)=\sum_{I\subset \{1,\dots,t\},I\ne
\emptyset}(-1)^{|I|+1}\mu(X\cap(\cap_{i\in I}B_{\nu_i})).
\end{eqnarray}
We have
\begin{eqnarray*}
X\cap Y=\cup_{l=1}^t(A_{\lam_1}\cap B_{\nu_l}),\\
X\cap (\cap_{i\in I}B_{\nu_i})=\cap_{i\in I}(A_{\lam_1}\cap
B_{\nu_i}).
\end{eqnarray*}
If $A_{\lam_1}\ne B_{\nu_l}$ for any $l$, then the type of
$A_{\lam_1}\cap B_{\nu_l}$ is strictly less than $r$, and
(\ref{add}) follows by the additivity of $\mu$ on $\ct_{r-1}$.
Assume now that $A_{\lam_1}=B_{\nu_1}$. Then $\mu(X\cap
Y)=\mu(B_{\nu_1})$. Also the right hand side in (\ref{add}) is
equal to
\begin{eqnarray*}
\mu(B_{\nu_1})+\sum_{I\subset \{2,\dots,t\},I\ne
\emptyset}(-1)^{|I|+1}\left(\mu(X\cap (\cap_{i\in I}B_{\nu_i}))-
\mu(X\cap B_{\nu_1}\cap(\cap_{i\in
I}B_{\nu_i}))\right)=\mu(B_{\nu_1}).
\end{eqnarray*}
This proves (\ref{add}).

Next we have
$$\mu(X\cup Y)=\mu(A_{\lam_1}\cup(\cup
_{l=1}^tB_{\nu_l}))=$$ $$\sum_{I\subset\{1,\dots
,t\},I\ne\emptyset}(-1)^{|I|+1}\mu(\cap_{i\in I}
B_{\nu_i})+\mu(A_{\lam_1})+\sum_{I\subset\{1,\dots
,t\},I\ne\emptyset}(-1)^{|I|+2}\mu(A_{\lam_1}\cap (\cap_{i\in I}
B_{\nu_i}))=$$
$$
\mu(Y)+\mu(X)-\mu(X\cap Y).$$

Let us assume that $s>1$. Then let us present $X=F\cup G$ where
$F$ and $G$ can be presented as a union of a smaller number than
$s$ of elements of $\cd$. Then by the assumption of induction we
have
$$\mu(X\cup Y)=\mu(F\cup (G\cup Y))=\mu(F)+\mu(G\cup Y)-\mu(F\cap(G\cup Y))=$$
$$\mu(F)+( \mu(G)+\mu(Y)-\mu(G\cap Y))-(\mu(F\cap G)+\mu(F\cap
Y)-\mu(F\cap G\cap Y))=$$ $$\mu(F\cup G)+\mu(Y)-\mu((F\cup G)\cap
Y) =\mu(X)+\mu(Y)-\mu(X\cap Y).$$ Thus $\mu$ is indeed a finitely
additive measure. \qed
\subsection{The sheaf  of finitely additive
measures.}\label{sheaf1} Let $X$ be a smooth manifold (of class
$C^{\infty}$). Let $\cp(X)$ denote the family of all simple
regular subpolyhedra of $X$ in sense of Definitions \ref{poly-1}
and \ref{poly-7}(1). Let $n=\dim X$.
\begin{definition}\label{sheaf-def}
A finitely additive measure $\mu$ on $\cp(X)$ is a functional
$$\mu:\cp(X)\to \CC$$
which satisfies the following property. Fix any $P\in \cp(X)$ and
any simple subdivision $\cd=\{P_{\lam}\}_{\lam\in \Lam}$ of $P$.
Define a function $m:\Lam\to \CC$ by $m(\lam):=\mu(P_\lam)$. For
$r=0,\dots,r$ let $\Lam_r:=\{\lam\in \Lam|\, \dim P_\lam =r\}$.
Then clearly $\Lam=\Lam_0\coprod \Lam_1\coprod \dots \coprod
\Lam_r$. The assumptions (1)-(3) before Lemma 2.2.2 are satisfied.
Let $\ct$ denote the family of all subsets representable as finite
unions of elements of $\cd$. Clearly $\ct$ is a finite family
closed under (finite) unions and intersections, and $P\in \ct$.
Let $\mu'$ denote the finitely additive measure on $\ct$ which is
constructed from $m$ as in Lemma \ref{fa-2}. Then we call $\mu$ to
be a finitely additive measure on $\cp(X)$ if $\mu(P)=\mu'(P)$ for
any $P$ and any subdivision $\cd$ of it.
\end{definition}

The linear space of all finitely additive measures on $\cp(X)$ we
will denote by $\cm(X)$. Now let us consider a presheaf $\cm_X$ of
vector spaces on $X$ defined as follows. For any open subset
$U\subset X$ set
$$\cm_X(U):=\cm(U)$$
with the obvious maps of restriction.
\begin{proposition}\label{sheaf-prop}
The presheaf $\cm_X$ is a sheaf.
\end{proposition}
{\bf Proof.} Let $U$ be any open subset of $X$. Let $\{U_\alp\}$
be any open covering of $U$. We have to check the following two
conditions:

(1) if $\mu\in \cm_X(U)$ is such that $\mu|_{U_\alp}=0$ for any
$\alp$ then $\mu=0$;

(2) if  we are given $\mu_\alp\in \cm_X(U_\alp)$ such that
$$\mu_\alp|_{U_\alp\cap U_\beta}=\mu_\beta|_{U_\alp\cap
U_\beta} \, \forall \alp,\,\beta$$ then there exists $\mu\in
\cm_X(U)$ such that $$\mu|_{U_\alp}=\mu_\alp \, \mbox{ for all }
\alp.$$

First let us check the condition (1). Let $P\in \cp(U)$. By
Proposition \ref{poly-17} we can choose a simple subdivision
$\cd=\{P_\lam\}$ of $P$ subordinate to $\{U_\alp\}$. Then one has
$$\mu(P)=\sum_{I\subset \Lam,\, I\ne\emptyset} (-1)^{|I|+1}
\mu(\cap_{i\in I}P_i)=0.$$

Let us check condition (2). Let $P\in \cp(U)$. Let us choose any
subdivision $\cd=\{P_\lam\}$ of $P$ subordinate to the covering
$\{U_\alp\}$. Let us define a function $$m:\Lam\to \CC$$ as
follows.  Let $\lam\in \Lam$. Choose $U_\alp$ such that
$P_\lam\subset U_\alp$. Define $m(\lam):=\mu_\alp(P_\lam)$.
Clearly $m$ is well defined. By Lemma \ref{fa-2} we can define a
number $\mu_\cd (P)$ using this subdivision $\cd$.
\begin{claim}\label{sheaf-claim}
The value $\mu_\cd (P)$ does not depend on the choice of a
subdivision $\cd$ of $P$.
\end{claim}
This value will be denoted by $\mu(P)$. Let us prove Claim
\ref{sheaf-claim}. Let $\cd_1$ and $\cd_2$ be two simple
subdivisions subordinate to the covering $\{U_\alp\}$. By
Proposition \ref{poly-18} we can choose simple subdivisions
$\cd_3,\, \cd',\, \cd''$ such that $\cd'$ is a refinement of
$\cd_1$ and $\cd_3$, and $\cd''$ is a refinement of $\cd_2$ and
$\cd_3$.

Thus in order to check that $\mu$ is well defined it remains to
check that if $\cd'$ is a refinement of $\cd$ then
$$\mu_\cd(P)=\mu_{\cd'}(P).$$
But this statement follows immediately from the uniqueness in
Lemma \ref{fa-2}.

To finish the proof of Proposition \ref{sheaf-prop} it remains to
prove that $\mu$ is indeed a finitely additive measure. Let $P\in
\cp(U)$. Let $\cd$ be any simple subdivision of $P$. Let $\cd'$ be
a simple refinement of $\cd$ subordinate to $\{U_\alp\}$. Then
define $m:\Lam\to \CC$ by $m(\lam)=\mu(P_\lam)$. The result
follows from Lemma \ref{fa-2}. \qed

\subsection{Smooth valuations.}\label{sheaf-sm} In this subsection we
introduce the main object of this article, namely smooth
valuations.

\def\an{CC}
Let $X$ be a smooth manifold of dimension $n$. Let $P\in \cp(X)$.
For any point $x\in P$ let us define the {\itshape tangent cone}
to $P$ at $x$, denoted by $T_xP$, the set
\begin{eqnarray}\label{tag-con}T_xP:=\{\xi\in T_xX| \mbox{ there exists a }
C^1-\mbox{map } \gamma\colon [0,1]\to P \mbox{ such that
}\gamma(0)=x\mbox{ and }\gamma'(0)=\xi\}.
\end{eqnarray}
It is easy to see that $T_xP$
coincides with the usual tangent space if $x$ is an interior point
of $P$. In general $T_xP\subset T_xX$ is a closed polyhedral cone.
\begin{definition}\label{define-char-c}
The {\itshape characteristic cycle} of $P$ is defined by
\begin{eqnarray}\label{cc-def}
CC(P):=\cup_{x\in P}(T_xP)^o
\end{eqnarray}
where for a convex cone $C$ in a linear space $W$ one denotes by
$C^o$ the dual cone
$$C^o:=\{y\in W^*|\, y(x)\geq 0\mbox{ for any } x\in C\}.$$
\end{definition}
Then $CC(P)$ is an $n$-dimensional subset of $T^*X$. It is
invariant under the the group $\RR_{>0}$ of positive real numbers
acting on $T^*X$ by the multiplication of cotangent vectors.
Moreover it is a Lagrangian submanifold with singularities. Note
that when $X$ is oriented the smooth part of $CC(P)$ carries an
induced orientation; then it is a cycle, i.e. $\pt (CC(P))=0$.
\begin{definition}\label{define-n-c}
The {\itshape normal cycle} $N(P)$ of $P$ is defined by
\begin{eqnarray}\label{nc-def}
N(P):=(a(CC(P))\backslash \{\underline{0}\})/\RR_{>0}
\end{eqnarray}
where $a\colon T^*X\to T^*X$ is the natural involution of
multiplication by $-1$ each cotangent vector, $\underline{0}$
denotes the zero section of $T^*X$, and the quotient is taken with
respect to the natural action of the group $\RR_{>0}$ mentioned
above. \end{definition}
Thus $N(P)\subset \PP_+(T^*X)$ is
$(n-1)$-dimensional submanifold with singularities. An orientation
of $X$ induces an orientation of $N(P)$; then it is a cycle. For
some references on the notions of the normal and characteristic
cycles see Remark \ref{remark-cc-nc}.

Let $\mu$ be a finitely additive measure on $X$ in sense of
Definition \ref{sheaf-def}.
\begin{definition}\label{ddd}
A measure $\mu$ is called {\itshape continuous valuation} if for
any  sequence of sets $\{P_N\}\subset \cp(X)$ which is contained
in a compact subset of $X$ and such that $\sup_N
vol(N(P_N))<\infty$, and a subset $P\in\cp (X)$ such that
$$\an(P_N)\to \an(P)$$
in sense of currents, one has $\mu(P_N)\to \mu(P)$.
\end{definition}
\begin{remark}
(1) Remind that the convergence in sense of currents means that
for any $\ome\in C^\infty(T^*X,\Ome^n\otimes p^*o)$ such that the
restriction of the projection $p$ to the support of $\ome$ is
proper, one has
$$\int_{CC(P_N)}\ome\to \int_{CC(P)}\ome.$$

(2) The convergence used in Definition \ref{ddd} is equivalent to
the flat convergence of currents, see \cite{federer-book}. The
equivalence is proved in \cite{simon}, Theorem 31.2.
\end{remark}
For any open subset $U\subset X$ let us denote by $\cc(U)$ the
space of continuous valuations on $U$. Clearly the correspondence
$U\mapsto \cc(U)$ is a sub-presheaf of $\cm_X$. It will be denoted
by $\cc_X$.

We would like to formulate a conjecture.
\begin{conjecture}
The presheaf $\cc_X$ is a sheaf.
\end{conjecture}

Let us denote by $\ck(\RR^n)$ the family of convex compact subsets
of $\RR^n$.

\begin{definition}
A measure $\mu$ is called {\itshape smooth valuation} if every
point $x\in X$ has a neighborhood $U\ni x$ and a diffeomorphism
$\phi:U\tilde\to \RR^n$ such that the restriction of $\phi_*\mu$
to $\cp(\RR^n)\cap \ck(\RR^n)$ extends by continuity in the
Hausdorff metric to $\ck(\RR^n)$ (clearly this extension is unique
if it exists) and this extension belongs to $SV(\RR^n)$ (see
Subsection \ref{valuations}).
\end{definition}

\def\vi{V^{\infty}}
\def\cvi{\cv^{\infty}}
For any open subset $U\subset X$ let us denote by $\vi(U)$ the set
of smooth valuations.

\begin{lemma}\label{sheaf-sm-van}
Let $V$ be an affine $n$-dimensional space. Let $\mu\in \cc(V)$.
Assume that $\mu(P)=0$ for any convex polytope $P$. Then $\mu=0$.
\end{lemma}
{\bf Proof.} Since by Proposition \ref{poly-8} any $P\in \cp(V)$
admits a triangulation, it is enough to show that $\mu$ vanishes
on any smoothly imbedded simplex $T$. Let $T=f(\Delta)$ where
$\Delta$ be the standard $n$-dimensional simplex in $\RR^n$, and
$f$ is a diffeomorphism of a neighborhood of $\Delta$ onto an open
subset in $V$. (The case of lower dimensional simplices in reduced
to $n$-dimensional case by approximation.) Let $K_N$ be a sequence
of convex compact subsets of $\RR^n$ with smooth boundary which
converges to $\Delta$ in the Hausdorff metric. Then $\an(K_N)\to
\an(\Delta)$, and $\sup_N vol(N(K_N))<\infty$ (this fact is known
for a long time, see e.g. the end of Section 1 in M. Z\"ahle
\cite{zahle90} where this fact was stated without proof; for a
proof we refer to \cite{part3} due to the lack of original
reference). Set $A_N:=f(K_N)$. Then $\an(A_N)\to \an (T)$ and
$\sup_N vol(N(A_N))<\infty$. Hence it is enough to show that
$\mu(A_N)=0$. Thus if one shows that for any compact domain $B$
with smooth boundary there exists a sequence of subsets $\{B_N\}$
presentable as a finite union of convex polytopes such that this
sequence is contained in a compact subset, $CC(B_N)$ have
uniformly bounded volume, and $CC(B_N)\to CC(B)$, then it follows
that $\mu(A_N)=0$. In this form this result is proved in
\cite{part3}; however the main step in the proof showing
convergence of the normal cycles (instead of the characteristic
cycles) is due to M. Z\"ahle \cite{zahle90}. The result follows.
\qed

\begin{lemma}\label{sheaf-sm-1}
Let $X$ be a smooth manifold.

(i) Let $\nu\in C^{\infty}(X,|\omega_X|)$, $\eta\in
C^{\infty}(\PP_+(T^*X),\Omega^{n-1}\otimes p^*o)$. Then $P\mapsto
\nu(P)+\int_{N(P)}\eta$ defines a smooth valuation on $X$.

(ii) Let $\mu\in \vi(X)$. Let $x\in X$. Then there exists a
neighborhood $U$ of $x$, $\nu\in C^{\infty}(U,|\omega_X|)$,
$\eta\in C^{\infty}(\PP_+(T^*U),\Omega^{n-1}\otimes o_{U})$ such
that for any $P\in \cp(U)$ one has
$$\mu(P)=\nu(P)+ \int_{N(P)}\eta.$$
\end{lemma}
{\bf Proof.} Part (i) follows from Proposition \ref{val-14}.

 Part (ii) immediately follows from Lemma
\ref{sheaf-sm-van} and Theorem \ref{val-15}. \qed

\begin{corollary}\label{sheaf-sm-2}
For any open subset $U\subset X$ the set of smooth valuations
$V^\infty(U)$ is a linear subspace of $\cm_X(U)$.
\end{corollary}
{\bf Proof.} This immediately follows from Lemma \ref{sheaf-sm-1}.
\qed


\begin{theorem}\label{sheaf-sm-3}
The correspondence $U\mapsto \vi(U)$ is a subsheaf of $\CC$-vector
spaces of the sheaf $\cm_X$.
\end{theorem}
This subsheaf will be denoted by $\cvi_X$.

{\bf Proof of Theorem \ref{sheaf-sm-3}.} By Corollary
\ref{sheaf-sm-2} $V^\infty(U)$ is a $\CC$-linear subspace of
$\cm_X(U)$. It immediately follows from Lemma \ref{sheaf-sm-1}
that $\cvi_X$ is a presheaf. Since the definition of $V^\infty(U)$
is local, the sheaf property of $\cvi_X$ is satisfied
automatically. \qed

Further properties of smooth valuations will be studied in the
next section. Now we will remind the following well known lemma
(see p.234 in \cite{schneider-book}; compare with Theorem 1.8.8 of
\cite{schneider-book}).
\begin{lemma}\label{4.8}
Let us fix a Euclidean metric on an affine space $V$. Let
$\{K_N\}$ be a sequence in $\ck(V)$ converging in the Hausdorff
metric to $K\in \ck(V)$. Let $A\in \ck(V)$. Then for almost all
isometries $g$ of $V$ one has
$$K_N\cap (gA)\to K\cap (gA)$$
in the Hausdorff metric.
\end{lemma}

\begin{proposition}\label{4.9}
Let $V$ be a linear space.

(i) The restriction map
$$\cc(V)\to CV(V)$$ is injective.

(ii) Under the above imbedding the image of $V^\infty(V)$ is equal
to $SV(V)$.
\end{proposition}
{\bf Proof.} Part (i) follows immediately from Lemma
\ref{sheaf-sm-van}.

Let us prove part (ii). First observe that Theorem \ref{val-15}
and Lemma \ref{sheaf-sm-1}(1) imply immediately that $SV(V)$ is
contained in the image of $V^\infty(V)$. To prove the opposite
inclusion let us fix a Euclidean metric on $V$ and fix $\phi\in
V^\infty (V)$. Let $\{U_\alp\}$ be an open covering of $V$ such
that, as in Lemma \ref{sheaf-sm-1}(ii), for any $\alp$ and any
$P\subset U_\alp$ $$\phi(P)=\nu_\alp(P)+\int_{N(P)}\eta_\alp$$
where $\nu_\alp\in C^\infty(U_\alp,|\ome_\alp|),\, \eta_\alp\in
C^\infty(\PP_+(T^*U_\alp),\Ome^{n-1}\otimes p^*o)$.

Let $K_0\in \ck(V),\, t_0\in [0,1],\, x_0\in V$. Assume first that
there exists $\alp_0$ such that $K_0+t_0+x_0\subset U_{\alp_0}$.
Then there exist neighborhoods $\co_1\subset [0,1]$ of $t_0$ and
$\co_2\subset V$ of $x_0$ such that for any $t\in \co_1$ and any
$x\in \co_2$ one has
$$K_0+t+x\subset U_{\alp_0}.$$ Then the function $[(t,x)\mapsto
\phi(tK_0+x)]$ is infinitely smooth in $\co_1\times \co_2$ by
Proposition \ref{val-14}(i).

Now an arbitrary $K_0\in \ck(V)$ can be represented as a finite
union of convex compact sets $K_0=\cup_{l=1}^sK_l$ such that for
each $l=1,\dots,s$ the set $K_l+t_0+x_0$ is contained in some
element of the covering $\{U_\alp\}$. The inclusion-exclusion
property and the above case imply the smoothness of the function
$[(t,x)\mapsto \phi(tK+x)]$ where $(t,x)\in [0,1]\times V$.

Thus it remains to show that the map $\ck(V)\to
C^\infty([0,1]\times V)$ given by $K\mapsto [(t,x)\mapsto
\phi(tK+x)]$ is continuous. Let us fix a lattice $L\subset V$. Let
$Q$ be a unit parallelepiped for $L$. It is easy to see that if
$\eps>0$ is small enough, then for any $x\in \eps L$ the set
$K\cap (x+\eps Q)$ is contained in one of the elements of the
covering $\{U_\alp\}$. Also we have $K=\cup_{x\in \eps
L}(K\cap(x+\eps Q))$. Replacing $L$ by its image under generic
isometry of $V$ close to the identity and using Lemma \ref{4.8} we
may assume that for  any $x\in \eps L$
$$K_N\cap (x+\eps Q)\to K\cap (x+\eps Q)$$
and similar convergence holds for finite intersections of the
above sets. Now the result follows from Proposition
\ref{val-14}(i). \qed

\section{Further properties of smooth valuations.}\label{further}
In Subsection \ref{further-1} we introduce and study the
filtration $\cw_\bullet$ on smooth valuations; in Proposition
\ref{filt-9} we obtain a description of smooth valuations in terms
of the integration with respect to the characteristic cycle. In
Subsection \ref{topology} we define the natural structure of
Fr\'echet space on the space of smooth valuations. In Subsection
\ref{eu-ve} we introduce the Euler-Verdier involution on smooth
valuations. \setcounter{subsection}{0}
\subsection{Filtration on smooth valuations.}\label{further-1}
\begin{definition}\label{filt-5}
Let $0\leq i\leq n$. Let $U$ be an open subset of a manifold $X$.
Let us denote by $\cw_i(U)$ the subset of $\cvi_X(U)$ consisting
of all elements $\phi\in \cvi_X(U)$ such that for every point
$x\in U$ there exists a neighborhood $V$ and a diffeomorphism
$f:V\tilde \to \RR^n$ such that the image of $f_*\phi$ in
$SV(\RR^n)$ belongs to $W_i$ (see Subsection \ref{valuations}).
\end{definition}
\begin{proposition}\label{filt-6}
(i) For any $0\leq i\leq n$ and for any open subset $U\subset X$,
$\cw_i(U)$ is a vector subspace of $\cvi_X(U)$.

(ii) The correspondence $U\mapsto \cw_i(U)$ is a subsheaf of
$\cvi_X$.
\end{proposition}
{\bf Proof.} Let us fix an open subset $V\subset U$ and a
diffeomorphism $f:U\tilde \to \RR^n$. By Theorem \ref{val-20} and
Lemma \ref{sheaf-sm-van} any valuation $\phi$ on $V$ such that
$f_*\phi$ lies in $W_i(V)$ has the following form: there exists
$\eta\in \tilde C^\infty(T^*U,W_i(T^*U)\otimes p^*o)$ such that
$$\phi(P)= \int _{CC(P)}\eta.$$
Obviously the set of valuations having the above form is a vector
subspace of $\cvi_X(V)$.
This proves part (i) of the proposition.

The same reasoning implies that $U\mapsto \cw_i(U)$ is a
sub-presheaf of ${\cal C}_X$. Since the definition of $\cw_i$ is
local, it is a sheaf. \qed

Remind that by Proposition \ref{4.9} we have the identification
$\cv^\infty_V(V)=SV(V)$. Using this identification we have the
following proposition.

\begin{proposition}\label{filt-6.5}
Let $V$ be a linear space. Then $\cw_i(V)=W_i$.
\end{proposition}
{\bf Proof.} It is clear from the definition that $W_i\subset
\cw_i(V)$. Let us prove that $\cw_i(V)\subset W_i$. Let $\phi\in
\cw_i(V)$. Fix $K\in \ck(V),\, x\in V$. There exists a
neighborhood $\cu$ of $x$ and a diffeomorphism $f\colon \cu\tilde
\to \RR^n$ such that $f_*\phi\in W_i(\RR^n)$. By Theorem
\ref{val-20} there exists $\eta\in \tilde C^\infty
(T^*\RR^n,W_i(\Omega^n\otimes p^*o))$ such that for any $A\in
\ck(\RR^n)$ one has
\begin{eqnarray}\label{eto}
(f_*\phi)(A)=\int_{CC(A)}\eta.
\end{eqnarray}
By Proposition \ref{4.9}(i) the formula (\ref{eto}) still holds
for any $A\in \cp(\RR^n)$. Hence for $0\leq t \ll 1$ one has
$$\phi(tK+x)=\int_{CC(tK+x)}\bar f^*\eta$$
where $\bar f$ is the natural lift of $f$ to $T^*\cu$. Set $\omega
:=\bar f^*\eta\in \tilde C^\infty(T^*\cu, W_i(\Omega^n\otimes
p^*o))$. Theorem \ref{val-20} implies that
$$\int_{CC(tK+x)}\omega= O(t^i).$$
Hence $\phi\in W_i$. \qed



Let us introduce more notation. Let us consider the following
sheaf on $X$:
\begin{eqnarray}\label{7.5}
\cw'_i(U)&:=& \tilde C^\infty(T^*U,W_i(\Ome^n\otimes p^*o)),\,
i=0,\dots, n.
\end{eqnarray}
Integration with respect to the normal cycle defines the following
morphism of sheaves which, by Theorem \ref{val-20} and Proposition
\ref{filt-6.5}, is an epimorphism:
$$\cw'_i \twoheadrightarrow\cw_i.$$
Clearly $\cw'_i/\cw'_{i+1}$ is isomorphic to the sheaf $[U\mapsto
\tilde C^\infty(U,W_i(\Ome^n\otimes p^*o)/W_{i+1}(\Ome^n\otimes
p^*o))]$. Hence we have a continuous epimorphism
\begin{eqnarray}\label{d-xi}
\Xi_i\colon \tilde C^\infty(T^*U,W_i(\Ome^n\otimes
p^*o)/W_{i+1}(\Ome^n\otimes p^*o))\twoheadrightarrow
\cw_i(U)/\cw_{i+1}(U).
\end{eqnarray}

Next we have a continuous map
\begin{eqnarray}\label{d-psi}
\Psi_i\colon \tilde C^\infty (T^*U,\Omega^{n-i}_{T^*U/U}\otimes
p^*(\wedge^iT^*U)\otimes p^*o)\to C^\infty(U,Val_i(TU)).
\end{eqnarray}
This map $\Psi_i$ is defined pointwise
$$\Psi_i|_x\colon \tilde
C^\infty(T^*_xX,\Omega^{n-i}_{T^*_xX}\otimes \wedge^iT^*_xX\otimes
p^*o)\to Val_i^{sm}(T_xX)$$ using the integration with respect to
the characteristic cycle of a subset of $T_x^*X$.

By Lemma \ref{val-19} there exists a canonical isomorphism of
vector bundles
$$J_i\colon W_i(\Ome^n\otimes p^*o)/W_{i+1}(\Ome^n\otimes p^*o)\tilde\to
\Omega^{n-i}_{T^*U/U}(T^*U)\otimes p^*(\wedge^iT^*U)\otimes
p^*o.$$

Let us denote for brevity by $\underline{Val}_i(TX)$ the sheaf
$[U\mapsto C^\infty(U,Val_i(TU))]$. Define sheaves $\car_i, \cs_i$
by
\begin{eqnarray*}
\car_i(U):= \tilde C^\infty(T^*U,W_i(\Ome^n\otimes p^*o)/W_{i+1}(\Ome^n\otimes p^*o)),\\
\cs_i(U):= \tilde C^\infty (T^*U,\Omega^{n-i}_{T^*U/U}\otimes
p^*(\wedge^iT^*U)\otimes p^*o).
\end{eqnarray*}

We have the canonical map
$$\cw_i(U)/\cw_{i+1}(U)\to (\cw_i/\cw_{i+1})(U).$$ The composition of
this map with $\Xi_i$ from (\ref{d-xi}) gives a map
\begin{eqnarray}\label{d-txi}
\tilde \Xi_i(U)\colon \car_i(U)\to (\cw_i/\cw_{i+1})(U).
\end{eqnarray}
This map is compatible with restrictions to open subsets. Hence we
obtain a morphism of sheaves
\begin{eqnarray}\label{d-ttxi}
\tilde \Xi_i\colon \car_i\to \cw_i/\cw_{i+1}.
\end{eqnarray}

\begin{lemma}\label{rw}
(i) There exists a natural isomorphism of sheaves
$$\cw_i'/\cw_{i+1}'\tilde \to \car_i.$$

(ii) For any open subset $\cu\subset X$
$$\car_i(\cu)=\cw_i'(\cu)/\cw_{i+1}'(\cu).$$
\end{lemma}
{\bf Proof.} Part (i) is obvious. To prove part (ii) note that we
have an exact sequence of $\co_X$-modules
$$0\to \cw_i'\to \cw_{i+1}'\to \car_i\to 0.$$
Hence from the long exact sequence we get
$$0\to \cw_i'(\cu)\to \cw_{i+1}'(\cu)\to \car_i(\cu)\to
H^1(\cu,\cw_i').$$ But since by Example \ref{sh-11}
$\co_X$-modules are acyclic we have $H^1(\cu,\cw_i')=0$. The
result follows. \qed
\begin{lemma}\label{epl}
The morphism $\tilde \Xi_i\colon \car_i\to \cw_i/\cw_{i+1}$ is an
epimorphism of sheaves.
\end{lemma}
{\bf Proof.} This follows immediately from the facts that
$\cw_i'\to \cw_i$ is an epimorphism, and $\car_i\simeq
\cw_i'/\cw_{i+1}'$. \qed

 We will need the following proposition.
\begin{proposition}\label{diagr}
(i) There exists unique morphism of sheaves on $X$
$$I_i\colon \cw_i/\cw_{i+1}\to \underline{Val}_i(TX)$$
which makes the following diagram commutative:
\def\aa{\car_i}
\def\bb{\cw_i/\cw_{i+1}}
\def\cc{\cs_i}
\def\dd{\underline{Val}_i(TX)}
\def\ff{\tilde\Xi_i}
\def\gg{J_i}
\def\hh{I_i}
\def\kk{\Psi_i}
\begin{eqnarray}\label{filt-4.5}
\square<1`1`1`1;1200`700>[\aa`\bb`\cc`\dd;\ff`\gg`\hh`\kk]
\end{eqnarray}
(ii) This morphism $I_i$ is an isomorphism of sheaves.
\end{proposition}
{\bf Proof.} The uniqueness of the morphism $I_i$ follows
immediately from the surjectivity of $\tilde\Xi_i$.

Let us prove the existence. Observe first of all that for any open
subset $\cu\subset \RR^n$ using Lemma \ref{rw} we have
$$\car_i(\cu)=\cw_i'(\cu)/\cw_{i+1}'(\cu)\overset{\Xi_i}{\to}
\cw _i(\cu)/\cw_{i+1}(\cu).$$ This shows that
\begin{eqnarray}\label{oo}
\Xi_i(\car_i(\cu))\inj \cw _i(\cu)/\cw_{i+1}(\cu)\to
(\cw_i/\cw_{i+1})(\cu).
\end{eqnarray}

Let us fix an open subset $\cu\subset X$ diffeomorphic to $\RR^n$,
and fix a diffeomorphism $f\colon \cu\tilde \to \RR^n$. By
Proposition \ref{filt-6.5} $\cw_j(\RR^n)=W_j\subset SV(\RR^n)$.
Theorem \ref{val-20} and (\ref{oo}) imply that
$\Xi_i(\car_i(\RR^n))=W_i/W_{i+1}$. Let us construct a map denoted
also $I_i\colon W_i/W_{i+1}\to Val_i^{sm}(\RR^n)$ which makes the
following diagram commutative:
\def\ff{\Xi_i}
\def\gg{J_i}
\def\hh{I_i}
\def\kk{\Psi_i}
\begin{eqnarray}\label{diag2}
\square<1`1`1`1;1200`700>[\car_i(\RR^n)` W_i/W_{i+1}`
\cs_i(\RR^n)`
C^\infty(\RR^n,Val_i(T\RR^n))=C^\infty(\RR^n,Val_i^{sm}(\RR^n));
\ff`\gg`\hh` \kk]
\end{eqnarray}
As in Subsection \ref{valuations} for $\phi\in W_i/W_{i+1}$ let us
define
$$(I_i \phi)(x,K)=\lim_{r\to +0}\frac{1}{r^i}\phi(rK+x)$$
where $K\in \ck(\RR^n),\, x\in \RR^n$. Let us show that the
diagram (\ref{diag2}) is commutative. For $\phi\in W_i/W_{i+1}$
and for all $x\in V,\, K\in \ck(V)$  one has
$$(I_i\phi)(x,K)=\lim_{r\to +0}\frac{1}{r^i}\phi(rK+x).$$
Let us fix $\eta\in \tilde C^\infty(T^*\RR^n,W_i(\Ome^n\otimes
p^*o )/W_{i+1}(\Ome^n\otimes p^*o))$. Let us fix a basis
$e_1^*,\dots, e_n^*$ in $V^*$. Then we can write
$$J_i(\eta)=\sum_{j_1,\dots,j_i}\eta_{j_1,\dots,j_i}\otimes
e_{j_1}^*\wedge \dots\wedge e_{j_i}^*$$ where
$\eta_{j_1,\dots,j_i}\in
C^\infty(T^*\RR^n,\Omega^{n-i}_{T^*\RR^n/X}\otimes p^*o)$. Then
\begin{eqnarray*}
(I_i(\Xi_i\eta))(x,K)&=&\sum_{j_1,\dots,j_i}\lim_{r\to
+0}\frac{1}{r^i} \int_{CC(rK+r)}\eta_{j_1,\dots,j_i}\otimes
e_{j_1}^*\wedge \dots\wedge e_{j_i}^*\\
&=&\sum_{j_1,\dots,j_i}\int_{CC(K)}\eta_{j_1,\dots,j_i}|_{p^{-1}(x)}\otimes
e_{j_1}^*\wedge \dots\wedge e_{j_i}^*\\
&=&(\Psi_i(J_i\eta))(x,K).
\end{eqnarray*}

Thus the diagram (\ref{diag2}) is commutative. Pulling the diagram
(\ref{diag2}) back to $\cu$ and using Proposition \ref{filt-6.5}
we obtain commutative diagram of vector spaces
\def\ff{\Xi_i}
\def\gg{J_i}
\def\hh{I_i}
\def\kk{\Psi_i}
\begin{eqnarray}\label{diag3}
\square<1`1`1`1;1200`700>[\car_i(\cu)` \cw_i(\cu)/\cw_{i+1}(\cu)`
\cs_i(\cu)` C^\infty(\cu,Val_i(T\cu)); \ff`\gg`\hh` \kk]
\end{eqnarray}
Note however that the map $I_i$ in (\ref{diag3}) might depend on a
choice of a diffeomorphism $f\colon \cu\tilde \to \RR^n$. This
however does not happen due to the uniqueness of $I_i$ which has
been proved.

Thus we have constructed, by now, for every open subset
$\cu\subset X$ diffeomorphic to $\RR^n$ the unique map $I_i\colon
\cw_i(\cu)/\cw_{i+1}(\cu)\to C^\infty(\cu,Val_i(T\cu))$ which
makes the diagram (\ref{diag3}) commutative. Since
$\cw_i/\cw_{i+1}$ is the sheafification of the presheaf $[U\mapsto
\cw_i(U)/\cw_{i+1}(U)]$ this defines in unique way the map of
sheaves $I_i\colon \cw_i/\cw_{i+1}\to \underline{Val}_i(TX)$ which
makes the diagram (\ref{filt-4.5}) commutative. This proves part
(i) of the proposition.

Part (ii) follows immediately from Theorem \ref{val-18} and the
description of $I_i$ after it. \qed

 We would like to state separately the following immediate
 corollary.
\begin{corollary}\label{filt-7}
The quotient sheaf $\cw_i/\cw_{i+1}$ is canonically isomorphic to
the sheaf $\underline{Val}_i(TX)$.
\end{corollary}



\begin{proposition}\label{filt-8}
The sheaves $\cw_i$ are soft. In particular the sheaf $\cvi_X$ is
soft.
\end{proposition}
{\bf Proof.} Consider the filtration of $\cw_i$ by subsheaves
$$\cw_i\supset \cw_{i+1}\supset \dots\supset \cw_n.$$
By Corollary \ref{filt-7} $\cw_j/\cw_{j+1}$ is an $\co_X$-module
for any $j$. Hence by Example \ref{sh-11} $\cw_j/\cw_{j+1}$ is
soft. Hence $\cw_i$ is also soft by Proposition \ref{sh-13}. \qed

\begin{proposition}\label{filt-9}
For any section $\phi\in \Gamma(X,\cw_i)$ there exists $\eta\in
\tilde C^\infty(T^*X,W_i(\Ome^n\otimes p^*o))$ such that for any
$P\in \cp(X)$
$$\phi(P)=\int_{CC(P)}\eta.$$
\end{proposition}
{\bf Proof.} Remind that in (\ref{7.5}) we have introduced the
sheaves $\cw'_i$. We have the canonical epimorphism
$$\cw'_i\surj \cw_i.$$ We have to show that the map
$$\Gamma(X,\cw'_i)\to \Gamma(X,\cw_i)$$
is an epimorphism. For $i=n$ this is obvious.

Since the sheaves $\cw'_{j+1}$ and $\cw_{j+1}$ are soft we have
\begin{eqnarray*}
\Gamma(X,\cw_j'/\cw_{j+1}')=\Gamma(X,\cw_j')/\Gamma(X,\cw_{j+1}'),\\
\Gamma(X,\cw_j/\cw_{j+1})=\Gamma(X,\cw_j)/\Gamma(X,\cw_{j+1}).
\end{eqnarray*}
By the descending induction in $i$ it is enough to show that the
induced maps
$$\Gamma(X,\cw'_j/\cw'_{j+1})\to \Gamma(X,\cw_j/\cw_{j+1})$$
are epimorphisms for all $j$. We may assume that $j<n$. We have
seen that the morphism
$$\tilde\Xi_i\colon \cw'_j/\cw'_{j+1} \to \cw_j/\cw_{j+1}$$
is an epimorphism of sheaves. Moreover the sheaf
$\cw'_j/\cw'_{j+1}$ is an $\co_X$-module. But
$\tilde\Xi_i=I_i^{-1}\circ \Psi_i\circ J_i$, and $\Psi_i$ and
$J_i$ are morphisms of $\co_X$-modules. Set $\ck:=Ker
\tilde\Xi_i$. Hence $\ck$ is isomorphic (via $I_i^{-1}$) to an
$\co_X$-module. Hence by Example \ref{sh-11} $H^i(X,\ck)=0$ for
$i>0$. From the long exact sequence we have
$$\Gamma(X,\cw'_j/\cw'_{j+1})\to \Gamma(X,\cw_j/\cw_{j+1})\to
H^1(X,\ck)=0.$$ Thus Proposition \ref{filt-9} follows. \qed
\begin{corollary}\label{filt-10}
For any $\phi\in \cw_n(X)$ there exists $\nu\in
C^\infty(X,|\ome_X|)$ such that for any $P\in \cp(X)$
$$\phi(P)=\nu(P).$$
Moreover for any $i=0,1,\dots,n-1$ and any $\phi\in \cw_i(X)$
there exist $\nu\in C^\infty(X,|\ome_X|)$ and $\ome\in
C^\infty(\PP_+(T^*X), W_i(\Ome^{n-1})\otimes p^*o)$ such that for
any $P\in \cp(X)$
$$\phi(P)=\nu(P)+\int_{N(P)}\ome.$$
\end{corollary}

\subsection{Linear topology on smooth valuations.}\label{topology}
Let us describe the canonical Fr\'echet space structure on the
space of smooth valuations. By Corollary \ref{filt-10} we have an
epimorphism of linear spaces
$$\Theta\colon C^\infty(X,|\omega_X|)\bigoplus
C^\infty(\PP_+(T^*X),\Omega^{n-1}\otimes p^*o) \twoheadrightarrow
\cv^\infty_X(X).$$ The source space has a canonical Fr\'echet
space structure.  It is easy to see that the kernel of $\Theta$ is
closed. Let us define the topology on $\cv^\infty_X(X)$ as the
quotient topology. This is a Fr\'echet topology. By the same
argument we define a Fr\'echet topology on $\cv^\infty_X(U)$ for
any open subset $U\subset X$. The following proposition is
trivial.
\begin{proposition}\label{topology-1}
For any open subsets $U\subset V\subset X$ the restriction map
$\cv^\infty_X(V)\to\cv^\infty_X(U)$ is continuous.
\end{proposition}
\begin{proposition}\label{topology-2}
Let $V$ be an $n$-dimensional linear space. Consider the
isomorphism of linear spaces $\cv^\infty_V(V)\tilde \to SV(V)$
from Proposition \ref{4.9}(ii).

Then this is an isomorphism of Fr\'echet spaces.
\end{proposition}
{\bf Proof.} By the Banach inversion theorem it is enough to check
that the map $\cv^\infty_V(V) \to SV(V)$ is continuous. This is
clear from the definitions. \qed

\begin{proposition}\label{closed}
For any $i=0,1,\dots,n$, $\cw_i(X)$ is a closed subspace of
$\cv^\infty_X(X)$.
\end{proposition}
{\bf Proof.} The definition of $\cw_i$ and Corollary \ref{filt-10}
imply that a smooth valuation $\phi\in \cv_X^\infty(X)$ belongs to
$\cw_i(X)$ if and only if for any open subset $U\subset X$
diffeomorphic to $\RR^n$, any diffeomorphism $f\colon U\tilde\to
\RR^n$, any $K\in \ck(\RR^n)$, and any $x\in \RR^n$ one has
$$\frac{d^k}{dt^k}\big|_{t=0}(f_*\phi)(tK+x)=0 \mbox{ for } k<i.$$
It is easy to see that $\phi\mapsto
\frac{d^k}{dt^k}\big|_{t=0}(f_*\phi)(tK+x)$ is a continuous linear
functional on $\cv^\infty_X(X)$ for any $U,f,K,x,k$ as above.
Hence $\cw_i(X)$ is a closed subspace of $\cv^\infty_X(X)$. \qed

\subsection{The Euler-Verdier involution.}\label{eu-ve} In this subsection we
construct a canonical continuous involution on the sheaf of smooth
valuations which we call the Euler-Verdier involution. Thus
$$\sigma: \cv_X^\infty\to \cv_X^\infty$$ satisfies $\sigma^2=Id$.
This involution preserves the filtration $\cw_\bullet$.

Let us describe the construction of $\sigma$. Remind that we have
the sheaf $\cw'_0$ on $X$ defined by
$$\cw'_0(U)=\tilde C^\infty(T^*U,\Omega^n\otimes p^*o)$$
where as previously the last space denotes the space of infinitely
smooth sections of the bundle $\Omega^n\otimes p^*o$ such that the
restriction of the projection $p$ to the support of these sections
is proper. By Proposition \ref{filt-9} we have epimorphism of
sheaves
$$\Theta:\cw'_0\to \cv_X^\infty.$$
On the space $T^*X$ we have the involution $a$ of multiplication
by -1 in each fiber of the projection $p:T^*X\to X$. It induces
involution $a^*$ of the sheaf $\cw'_0$.
\begin{proposition}
The involution $(-1)^n a^*$ factorizes (uniquely) to involution of
$\cv_X^\infty$ denoted by $\sigma$.
\end{proposition}
{\bf Proof.} We have to show that if $\omega\in\tilde
C^\infty(T^*U,\Omega^n\otimes p^*o)$ satisfies $\Theta(\ome)=0$
then $\Theta( a^*(\omega))=0$.

It is easy to see that for any $\omega\in \tilde
C^\infty(T^*X,\Omega^n\otimes p^*o)$ and any $P\in \cp(X)$ one has
\begin{eqnarray}\label{verdier-equal}
\int_{CC(P)}a^*\omega =(-1)^{n-\dim P}\left( \int_{CC(P)}\omega
-\int_{CC(\pt P)}\omega\right)
\end{eqnarray}
where $\pt P:=P\backslash int P$, and $int P$ if the relative
interior of $P$. The formula (\ref{verdier-equal}) immediately
implies the proposition. \qed

The following result is clear from the discussion above.
\begin{theorem}\label{eu-ve-th}
(i) The Euler-Verdier involution $\sigma$ preserves the filtration
$\cw_\bullet$.

(ii) The induced involution on $gr_\cw \cv_X^\infty \simeq
\underline{Val}_\bullet(TX)$ comes from the involution on the
bundle $Val(TX)$ defined as $\phi\mapsto [K\mapsto (-1)^{deg \phi}
\phi(-K)]$ for any $\phi\in Val(T_xX)$ for any $x\in X$, and where
$deg \phi$ is the degree of homogeneity of $\phi$.
\end{theorem}

Thus the sheaf $\cv_X^\infty$ of smooth valuations decomposes
under the action of the Euler-Verdier involution into two
subsheaves $\cv_X^{\infty, +}$ and $\cv_X^{\infty, -}$
corresponding to eigenvalues 1 and -1 of $\sigma$ respectively.
Thus
$$\cv_X^\infty=\cv_X^{\infty, +}\oplus \cv_X^{\infty, -}.$$

\end{document}